\documentclass[10pt]{amsart}
\usepackage{amsfonts}
\usepackage{color}
\usepackage[square,compress,comma, numbers, sort]{natbib}
\usepackage[colorlinks=true, citecolor=blue, linkcolor=blue]{hyperref}
\allowdisplaybreaks[4]
\usepackage{amssymb}
\usepackage{amsmath}

\definecolor{c20}{rgb}{0.,0.7,0.}
\definecolor{c30}{rgb}{0.,0.,1.}
\definecolor{c40}{rgb}{1,0.1,0.7}
\definecolor{c50}{rgb}{1,0,0}
\definecolor{c60}{rgb}{1,0.9,0.1}

\allowdisplaybreaks[4]

\def\cLa#1{\textcolor{black}{#1}}

\def\cl#1{\textcolor{c50}{#1}}

\def\cJ#1{\textcolor{black}{#1}}

\def\EHb#1{\textcolor{c30}{#1}}
\def\EHc#1{\textcolor{black}{#1}}
\def\kd#1{\textcolor{cyan}{#1}}
\def\kd#1{\textcolor{black}{#1}}

\newcommand{\abs}[1]{\left\lvert #1 \right\rvert}
\newcommand{\Abs}[1]{ \biggl \lvert #1 \biggr \rvert}

\newcommand{\E}[1]{\mathbb{E}\left\{#1\right\}}

\newcommand{\pk}[1]{\mathbb{P} \left \{#1 \right \} }

\newcommand{\R}{\mathbb{R}}

\newcommand{\N}{\mathbb{N}}
\newcommand{\inr}{\in \R}

\newcommand{\ldot}{,\ldots,}

\newcommand{\limit}[1]{\lim_{#1 \to   \infty}}

\newcommand{\BQN}{\begin{eqnarray}}
\newcommand{\EQN}{\end{eqnarray}}
\newcommand{\BQNY}{\begin{eqnarray*}}
\newcommand{\EQNY}{\end{eqnarray*}}
\newcommand{\BS}{\begin{sat}}
\newcommand{\ES}{\end{sat}}
\newcommand{\BT}{\begin{theo}}
\newcommand{\ET}{\end{theo}}
\newcommand{\BL}{\begin{lem}}
\newcommand{\EL}{\end{lem}}
\newcommand{\BK}{\begin{korr}}
\newcommand{\EK}{\end{korr}}

\newcommand{\BD}{\begin{de}}
\newcommand{\ED}{\end{de}}

\newcommand{\BIT}{\begin{itemize}}
\newcommand{\EIT}{\end{itemize}}
\newcommand{\BDI}{\begin{description}}
\newcommand{\EDI}{\end{description}}

\newcommand{\BRM}{\begin{remarks}}
\newcommand{\ERM}{\end{remarks}}

\newcommand{\BEL}{\begin{lem}}
\newcommand{\EEL}{\end{lem}}

\newtheorem{theo}{Theorem}[section]
\newtheorem{sat}[theo]{Proposition}
\newtheorem{de}[theo]{Definition}
\newtheorem{lem}[theo]{Lemma}

\newtheorem{example}[theo]{Example}
\newtheorem{korr}[theo]{Corollary}
\newtheorem{remark}[theo]{Remark}
\newtheorem{remarks}[theo]{Remarks}
\newtheorem{prop}[theo]{Proposition}

\newcommand{\nelem}[1]{{Lemma \ref{#1}}}
\newcommand{\neprop}[1]{{Proposition \ref{#1}}}
\newcommand{\netheo}[1]{{Theorem \ref{#1}}}

\newcommand{\prooftheo}[1]{ \textsc{\bf Proof of Theorem} \ref{#1}:}
\newcommand{\proofprop}[1]{\textsc{\bf Proof of Proposition} \ref{#1}:}
\newcommand{\prooflem}[1]{\textsc{\bf Proof of Lemma} \ref{#1}:}

\newcommand{\COM}[1]{}

\newcommand{\QED}{\hfill $\Box$}

\def\rw{\rightarrow}

\def\IF{\infty}

\def\LT{\left}
\def\RT{\right}
\def\H{\mathcal{H}}

\def\rw{\rightarrow}

\def\H{\mathcal{H}}

\def\vn{\varepsilon}
\def\Var{\text{Var}}

\def\II{\mathbb{I}}

\def\Oxt{\overline{X}(t)}
\def\Xgt{(X(t)+g(t))}

\def\TPi{\widetilde{\Pi}}
\def\Bu+#1{\mathcal{B}^{\varepsilon+}_{u}(#1)}

\def\Oxt{\overline{X}(t)}
\def\orh{\overleftarrow{\rho}}
\def\qa{Q_{t_0}}
\def\pa{w_{t_0}}

\begin{document}

\title{Extremes of  Threshold-Dependent Gaussian Processes}%[Extremes of Gaussian Processes]

\author{Long Bai}
\address{Long Bai, %School Mathematical Sciences and LPMC, Nankai University, Tianjin 300071, China, and
Department of Actuarial Science, %\\Faculty of Business and Economics\\
University of Lausanne
UNIL-Dorigny, 1015 Lausanne, Switzerland
}
\email{Long.Bai@unil.ch}

\author{Krzysztof D\c{e}bicki}
\address{Krzysztof D\c{e}bicki, Mathematical Institute,
	University of Wroc\l aw, pl. Grunwaldzki 2/4, 50-384 Wroc\l aw, Poland}
\email{Krzysztof.Debicki@math.uni.wroc.pl}

\author{Enkelejd  Hashorva}
\address{Enkelejd Hashorva, Department of Actuarial Science, University of Lausanne,\\
UNIL-Dorigny, 1015 Lausanne, Switzerland
}
\email{Enkelejd.Hashorva@unil.ch}

\author{Lanpeng Ji}
\address{Lanpeng Ji, Department of Actuarial Science, University of Lausanne\\
UNIL-Dorigny, 1015 Lausanne, Switzerland
}
\email{jilanpeng@126.com}

\bigskip

\date{\today}
 \maketitle

{\bf Abstract}: \EHc{In this contribution we are concerned with the asymptotic behaviour, as $u\to \IF$,
of $\pk{\sup_{t\in [0,T]} X_u(t)> u}$, where $X_u(t),t\in [0,T],u>0$ is a family of centered Gaussian processes with continuous trajectories.
A key application of our
\kd{findings concerns
$\pk{\sup_{t\in [0,T]} (X(t)+ g(t))> u}$, as $u\to\infty$,
}
for $X$ a centered Gaussian process and $g$ some measurable trend function.
Further applications include  the approximation of both the ruin time and the ruin probability of \cJ{the Brownian motion risk model with constant force of interest.}}

{\bf Key Words:}  Extremes; Gaussian processes;   fractional Brownian motion; ruin probability; ruin time.

{\bf AMS Classification:} Primary 60G15; secondary 60G70

\section{Introduction}
Let $X(t), t\ge 0 $ be a centered Gaussian process with  continuous trajectories. An important problem in applied and theoretical probability is the determination of the asymptotic behavior of
\BQN\label{Ab1}
p(u)=\pk{\sup_{t\in[0,T]}(X(t)+ g(t)) >u},\ \ \ u\to\IF
\EQN
for some $T>0$ and \cJ{$g(t), t\in[0,T]$ a bounded measurable} %trend
\kd{function.}
 For instance, if $g(t)= - ct$, then in the context of risk theory $p(u)$
\kd{has interpretation as}
the ruin probability over the finite-time horizon $[0,T]$. %, where $X(t)$ denotes the total aggregate loss up to time $t$ and $c>0$ is some constant premium rate,
\kd{Dually, in the context of queueing theory, $p(u)$ is related to the buffer overload problem;}
see e.g., \cite{DebRol02,  debicki2002ruin,DI2005,MR3091101,DEJ14}.

For the special case that $g(t)=0, t\in [0,T]$ the exact asymptotics of  \eqref{Ab1} is well-known for both locally stationary and general
non-stationary Gaussian processes, see e.g.,
\cite{PicandsB, Pit72, GennaBorell, Berman92, Pit96,AZI, hashorva2000extremes, MR2462286, debicki2011extremes, Pit20,ChengC,ChengA,Marek}.
Commonly, for  $X$ a centered  non-stationary Gaussian process it is assumed that the standard deviation function $\sigma$
is such that  \kd{$t_0=argmax_{t\in [0,T]} \sigma(t)$} is unique and $\sigma(t_0)=1$. Additionally, if  the correlation function $r$ and the standard deviation function $\sigma$ satisfy
 (hereafter $\sim$ means asymptotic equivalence)
\BQN\label{assump-corre}
1- r(s,t)\sim a\abs{t-s}^\alpha,\qquad 1- \sigma(\cJ{t_0}+t)\sim b\abs{t}^\beta %(1+ o(1))+o(\abs{t}^{\beta'})
,\qquad s,t\to t_0
\EQN
for some $a,b,\beta$ positive and $\alpha\in (0,2]$,
then we have (see \cite{Pit96}[Theorem D.3]) % or \cite{PP78}[Theorem 1]
\BQN \label{pu} p(u) \sim C_0 u^{(\frac 2 \alpha- \frac  2 \beta)_+} \pk{X(t_0)>u}, \quad u\to \IF,
\EQN
where $(x)_+=\max(0,x)$ and
\BQNY
 C_0 =
 \LT\{
\begin{array}{ll} a^{ {1}/{\alpha}}b^{-1/\beta}\Gamma(1/\beta+1)\mathcal{H}_{\alpha},&\hbox{if} \ \ \alpha<\beta,\\
\mathcal{P}_{\alpha}^{b/a},&  \hbox{if} \ \ \alpha=\beta,\\
1,&  \hbox{if} \ \ \alpha>\beta.
\end{array}
 \RT.
\EQNY
Here $\Gamma (\cdot)$ is the gamma function, and %$\mathcal{H}_\alpha$ and $\mathcal{P}_{\alpha }^{b/a}$
\BQNY
\mathcal{H}_\alpha= %\E{ \frac{\sup_{t \in \R } e^{W(t)}}{{\int_{t \in \R } e^{W(t)}dt}}} ,
\lim_{T\rightarrow\IF}\frac{1}{T} \E{\sup_{t\in [0,T]}e^{W(t)}},
\quad  \mathcal{P}_{\alpha }^{b/a} = \E{ \sup_{t \in [0,\IF) } e^{W(t)- b/a \abs{t}^\alpha}},
\ \text{with \ } W(t)=\sqrt{2} B_\alpha(t)- \abs{t}^\alpha,
\EQNY	
are the Pickands and Piterbarg constants, respectively, where $B_\alpha$ is  a standard fractional Brownian motion (fBm) with self-similarity index $\alpha/2\in (0,1]$, see \cite{shao1996bounds,Harper1, Harper2, debicki2016parisian, mi:17, LongL, SBK} for properties of both constants.

The more general case with  non-zero $g$  has also been considered in the literature; see, e.g., \cite{PP78, Deb99, PitSta01, DebRol02, EnkelejdJi2014Chi,MR2462285}. However, most of the aforementioned contributions treat only restrictive trend functions   $g$. For instance, in \cite{PP78}[Theorem 3]  a  H{\"o}lder-type condition for $g$  is assumed, which excludes important cases of $g$ that appear in applications. The restrictions are often so severe that simple cases such as the Brownian bridge with drift considered in Example \ref{exam1} below cannot be covered.\\ %\\
%The main goal of this paper is to discuss the asymptotics of $p_g(u)$ for general Gaussian processes $X$ and general trend fucntions $g$. For a large class of trend functions our results are completely new.\\
\kd{A} key difficulty when dealing with $p(u)$ is that $X+g$ is not \kd{a centered Gaussian process}.
It is however possible to get rid of the trend function $g$ since for any bounded  function $g$ and all  $u$ large
 \eqref{Ab1}   can be re-written as
\BQNY
p_T(u)= \pk{\sup_{t\in[0,T]}  X_u(t)>u }, \quad X_u(t)=\frac{X(t)}{1-g(t)/u}, \quad t\in[0,T].
\EQNY
Here $X_u$ is centered, however it depends on the threshold $u$, which complicates the analysis.\\
Extremes of threshold-dependent Gaussian processes $X_u(t),t\in \R$ have been already dealt with in several contributions, see e.g., \cite{debicki2002ruin,DI2005,HP2004, MR2462285,DHJ13a}.
%, \cLa{which, however, mainly focus on the infinite time horizon or just consider some restrictive cases over finite time horizon.}
 Our principal result in Theorem \ref{PreThm2} derives the asymptotics of $p_T(u)$ for quite general families of centered Gaussian processes $X_u$ under tractable assumptions on the variance and correlation functions of $X_u$. To this end,  in \netheo{PreThm1} we first derive the asymptotics of
 $$
 p_\Delta(u)=\pk{\sup_{t\in \Delta(u)}  X_u(t)>u }, \ \ \ u\to\IF
 $$
for some short compact intervals $\Delta(u)$.

Applications of our  main results include  derivation of Proposition  \ref{Thm3} for a class of  locally stationary Gaussian processes with trend and that of Proposition \ref{MainThm1} for a class of non-stationary Gaussian processes with trend, as well as those of their corollaries.
For instance,
a direct application of Proposition \ref{MainThm1} yields the asymptotics of \eqref{Ab1} for a non-stationary $X$  with standard deviation function $\sigma$ and correlation function $r$  satisfying \eqref{assump-corre} with $t_0=argmax_{t\in [0,T]} \sigma(t)$.
If further the trend function $g$ is continuous in a neighborhood of $t_0,$ $g(t_0)=\max_{t\in [0,T]} g(t)$ and
\BQN \label{eq:gtt0}
g(t) \sim g(t_0)-c|t-t_0|^\gamma,\quad t \rw t_0
\EQN
for some positive constants $c, \gamma$, then \eqref{pu} holds %for $p(u)=\pk{\sup_{t\in 0,T} (X(t)+g(t)) > u}$ where
with $C_0$  specified in \neprop{Thm2} and $\beta,u$  being substituted by $\min(\beta, 2\gamma)$
and $u-g(t_0)$  respectivelly.

\COM{
Further as shown in \cite{DebRol02}, in the applications since we just need to know the structure of the standard deviation function
$$\sigma_{X_u}(t)=\sqrt{Var(X_u(t)) }=\frac{\sigma (t)}{1-g(t)/u}, \quad
\sigma  (t)= \sqrt{Var(X(t)},$$
rather than the structure of $g(t)$,
 some more general cases where $g(t)$ may not satisfy \eqref{eq:gtt0} also can be included in \neprop{MainThm1}. An instance about Brownian motion risk model shows this senario in \neprop{exmodel}.\\
Moreover, when $X+g$ is a locally staionary Gaussian process, say, $\sigma(t)\equiv1$, and $g(t)$ is continuous with \eqref{eq:gtt0} satisfied, for the asymptotic behavior of \eqref{Ab1} $g(t)$ has a similarly contribution as $\sigma(t)$ in \eqref{pu}, which is showed in \neprop{Thm3}.\\
Some generalizations, such as the regularly varying dependence structure about $\sigma$ and $r$  in \cite{regularly2016},  also can be included in \netheo{PreThm2}.
}

Complementary, we investigate asymptotic properties of the first passage time (ruin time) of $X(t)+g(t)$ to $u$ on the finite-time interval $[0,T]$,
given the process has ever  exceeded $u$ during  $[0,T]$. \kd{In particular}, for
\BQN\label{tu}
\tau_u=\inf\{t\geq0:  X(t) > u-g(t)\},
\EQN
with $\inf\{\emptyset\}=\IF$, we are interested in the approximate distribution of $\tau_u|\tau_u\leq T$, as $u\rw\IF$.
Normal and exponential approximations of various Gaussian models have been discussed in \cite{MR2462285,HJ13,DHJ13a,DHJParisian,KEP2015}.
In this paper, we derive general results for the approximations of the conditional passage time in Propositions \ref{ruintimelocal}, \ref{ruintime}.
\cJ{The asymptotics of $p_\Delta(u)$ for  a short compact intervals $\Delta(u)$ displayed in \netheo{PreThm1} plays a key  role in the derivation of these results}. \\
%Another unexpected obtainment is generalization of the classical Piterbarg constants. Some interesting and useful properties are given at \neprop{prop1} and  {\bf Remark} \ref{pep}.\\
%Applications of our findings include the approximation of the first ruin time and the ruin probability of a Brownian motion risk model with constant force of interest over the infinite-time horizon.
Organisation of the rest of the paper: In Section 2, the tail asymptotics of the supremum of a family of centered Gaussian processes  \cJ{indexed by} $u$ are given. Several applications and examples are displayed in Section 3. \COM{Two applications and some examples are discussed in Section 4 and Section 5, respectively.} Finally, we present all the proofs in Section 4 and Section 5.

\section{Main Results} %Preliminaries
\COM{
As shown  in \cite{Pit96}, if $X(t), t\in[0,T]$ has non-constant variance $\sigma^2_X(t)$, then the maximum point of $\sigma^2_X(t)$ plays a key role in the asymptotics of \eqref{Ab1}. %and the local structure of it around the maximumpoint
Assume that $t=t_0\in(0,T)$ is the unique maximum point of $\sigma^2_X(t)$. Under some conditions on the local structure of the variance function and the correlation function of $X$ we have, for some interval $\Delta(u)$,
\BQNY
\pk{\sup_{t\in[0,T]}X(t)>u}\sim \pk{\sup_{t\in\Delta(u)}X(t_0+t)>u},\ \ \ u\to\IF.
\EQNY}
\COM{Since the trend function $g$ is bounded, %by the assumption that $g$ is bounded on $[0,T]$
we can re-write \eqref{Ab1} for all  $u$ large as
\BQNY
p_g(u)= \pk{\sup_{t\in[0,T]}  X_u(t)>u }, \quad X_u(t)=\frac{X(t)}{1-g(t)/u}, \quad t\in[0,T].
\EQNY
As shown in \cite{DebRol02}, the standard deviation function
$$\sigma_{X_u}(t)=\sqrt{Var(X_u(t)) }=\frac{\sigma (t)}{1-g(t)/u}, \quad
\sigma  (t)= \sqrt{Var(X(t)}$$
plays a key role in the asymptotics of $p_g(u)$.
 Assuming  that $m_u^*= \max_{t\in[0,T]}\sigma_{X_u}(t)>0$ for all large $u$, we have thus %. Then we have
\BQN \label{ptu}
\EHb{p( [0,T], X_u)}:=p_g(u)= \pk{\sup_{t\in[0,T]} \frac{X_u(t)}{m^*_u }>\frac{u}{m^*_u} },
\EQN
which is a convenient expression to deal with.\\
}
Let $X_u(t), t\in\cJ{\mathbb R},u>0$ be a family of threshold-dependent centered Gaussian processes with continuous trajectories,
 variance functions $\sigma^2_u$ and correlation functions $r_u$.
 \cJ{Our main results concern the asymptotics of \kd{slight generalization of $p_\Delta(u)$  and $p_T(u)$ for families of centered Gaussian processes $X_u$
satisfying some regularity conditions for variance and coavariance respectivelly.}}

%For some  closed interval $F\subset\R$ and $F\ni 0$, define function collection $C^*_0(F)$: \\
%If $f(\cdot) \in C^*_0(F)$, then $f(\cdot)$ is non-negative continuous over $\widetilde{F}$ where $\widetilde{F}\supset F$ is a single interval  on $\R$,
Let $C^*_0(E)$ be the set of continuous real-valued functions defined on the interval $E$ such that
$f(0)=0$ and for some $\epsilon_2>\epsilon_1>0$
\BQN\label{eq:fff}
\lim_{|t|\rw\IF,t\in E}f(t)/|t|^{\epsilon_1}=\IF,\ \ \ \lim_{|t|\rw\IF,t\in E}f(t)/|t|^{\epsilon_2}=0,
\EQN
if  $\sup\{x:x\in E\}=\IF $ or $\inf\{x:x\in E\}=-\IF$.\\
In the following $\mathcal{R}_\alpha$ denotes the  set of regularly varying functions at $0$ with index $\alpha\inr$, see \cite{EKM97,Res,Soulier} for details.\\
We shall impose the following assumptions where $\Delta(u)$ is a compact interval:\\
\textbf{A1}: For any large $u$, there exists a point $t_u\in \R$ such that $\sigma_u(t_u)=1$.\\
\textbf{A2}: There exists some
$\lambda>0$  such that
\begin{align}\label{A2}
\lim_{u\rightarrow\IF}\sup_{t \in\Delta(u)}\left|\frac{ \LT(\frac{1}{\sigma_u(t_u+t)}-1\RT)u^2-f(u^\lambda t)}{f(u^{\lambda}t)+1}\right|= 0
\end{align}
holds for some non-negative continuous function $f$ with $f(0)=0$.\\
\textbf{A3}: There exists $\rho\in\mathcal{R}_{\alpha/2}, \alpha\in(0,2]$ such that
$$
\lim_{u\rightarrow\IF}\underset{t\not=s }{\sup_{s,t\in\Delta(u)}}\left|\frac{1-r_u(t_u+s,t_u+t)}{\rho^2(|t-s|)}-1\right|=0
$$
\kd{and $
\eta:=\lim_{s\to 0} \frac{ \rho^2(s)}{ s^{2/\lambda}}\in[0,\IF]$,
with $\lambda$ given in  {\bf A2}.}

\begin{remark}\label{remarkVar}
If  $f$ satisfies $f(0)=0$ and $f(t)>0, t\neq 0$, then
\begin{align*}
\lim_{u\rightarrow\IF}\sup_{t \in\Delta(u), t\neq0}
\LT|\frac{\frac{1}{\sigma_u(t_u+t)}-1}
{u^{-2}f(u^{\lambda}t)}-1\RT|=0
\end{align*}
for some $\lambda>0$ implies that \eqref{A2} is valid.
\COM{In fact
\begin{align*}
0&\leq\lim_{u\rightarrow\IF}\sup_{t \in\Delta(u)}\left|\frac{ \LT(\frac{1}{\sigma_u(t_u+t)}-1\RT)u^2-f(u^\lambda t)}{f(u^{\lambda}t)+1}\right|= \lim_{u\rightarrow\IF}\sup_{t \in\Delta(u), t\neq0}\left|\frac{ \LT(\frac{1}{\sigma_u(t_u+t)}-1\RT)u^2-f(u^\lambda t)}{f(u^{\lambda}t)}\right|\left|\frac{ f(u^{\lambda}t)}{f(u^{\lambda}t)+1}\right|\\
&\leq \lim_{u\rightarrow\IF}\sup_{t \in\Delta(u), t\neq0}\left|\frac{ \frac{1}{\sigma_u(t_u+t)}-1}{u^{-2} f(u^{\lambda}t)}-1\right|=0.
\end{align*}
\begin{remark}
If non-negative continuous function $f(\frac{1}{t})\in\mathcal{R}_{-\epsilon}$ and $f(0)=0$, then $f\in C^*_0([0,\IF))$.
\end{remark}}
\end{remark}
\COM{It is apparent that the key in analysis of \eqref{Ab1} is to investigate
\BQN
\pk{\sup_{t\in\Delta(u)} X_u(t+t_u)>M_u},\ u\rw\IF,
\EQN
where $M_u$ satisfies $\lim_{u\rw\IF}\frac{M_u}{u}=c_0\in(0,\IF)$.}
Next  we  introduce some further notation, starting with the  Pickands-type constant   defined by
$$%\mathcal{H}_\alpha=\lim_{T\rightarrow\IF}\frac{1}{T}\mathcal{H}_\alpha[0,T],\quad \text{with }
\mathcal{H}_\alpha [0,T]=\E{\sup_{t\in [0,T]}e^{\sqrt{2}B_\alpha(t)-|t|^\alpha}},  \ \ T>0, %\in(0,\IF),
$$
where $B_\alpha $ is \kd{an} fBm.
Further, define for $f \in C^*_0([S,T])$ with $S,T\in\R, S<T$ and a positive constant $a$
\BQNY
&&\mathcal{P}_{\alpha,a}^f [S,T]=\E{\sup_{t\in [S,T]}e^{\sqrt{2a}B_\alpha(t)-a|t|^\alpha-f(t)}},
\EQNY
and set %for certain continuous function $f$ % $ f\in C^*_0(\R)$
\BQNY
&&\mathcal{P}_{\alpha,a}^f[0,\IF)=\lim_{T\rw\IF} \mathcal{P}_{\alpha,a}^f[0,T], \quad
\mathcal{P}_{\alpha,a}^f(-\IF,\IF)=\lim_{S\rw-\IF,T\rw\IF} \mathcal{P}_{\alpha,a}^f[S,T].
\EQNY
The finiteness of $\mathcal{P}_{\alpha,a}^f [0,\IF)$ and $\mathcal{P}_{\alpha,a}^f (-\IF,\IF)$ is guaranteed under weak assumptions on $f$, which will be shown in the proof of \netheo{PreThm1}, see \cite{PicandsA,Pit72, debicki2002ruin,DI2005,DE2014,DiekerY,DEJ14,Pit20, Tabis, DM, SBK} for various properties of $\mathcal{H}_{\alpha}$ and $\mathcal{P}_{\alpha,a}^f[0,\IF)$. \\
%In our notation, $\sim$ means asymptotic equivalence when the argument tends to $0$ (or $\IF$). Below $ \Phi(\cdot)$ and $\Psi(\cdot) $ stand for the distribution function and survival function of an $N(0,1)$ random variable, respectively. %; recall that $\Psi(u)\sim \frac{1}{\sqrt{2\pi}u}e^{-\frac{u^2}{2}}, u\rw\IF$.
Denote by %$\Gamma_{\{\cdot\}}$ the gamma function and
$\II_{\{\cdot\}}$ the indicator function.
%Denote by $\II_{\{\cdot\}}$ the indicator function.
For the regularly varying function $\rho(\cdot)$, we denote by $\overleftarrow{\rho}(\cdot)$ its asymptotic inverse (which is asymptotically unique).
Throughout this paper, we set $0\cdot\IF =0$ and  $u^{-\IF}=0$ if  $u>0$.
\kd{Let $\Psi(u):=\pk{\mathcal{N}>u}$, with $\mathcal{N}$ a standard normal random variable.}
\\
%Hereafter we  and
\COM{\BQN \lim_{u\rw\IF}f(u)[y_1(u),y_2(u)]:=[y_1,y_2],
\EQN
if $\lim_{u\rw\IF}f(u)y_1(u)=y_1\in[-\IF,\IF),\ \lim_{u\rw\IF}f(u)y_2(u)=y_2\in(-\IF,\IF]$ and $y_1<y_2$.\\}
In the next theorem  we shall consider two functions $x_1(u),x_2(u),u\inr$ such that $x_1(\frac{1}{t})\in\mathcal{R}_{\mu_1},\ x_2(\frac{1}{t})\in\mathcal{R}_{\mu_2}$ with $\mu_1,\mu_2\geq \lambda$, and
\BQN\label{x1x2}
\lim_{u\rw\IF}u^{\lambda}x_i(u)=x_i\in[-\IF,\IF],i=1,2, %\ \lim_{u\rw\IF}u^{\lambda}x_2(u)=x_2\in(-\IF,\IF],\
\quad \hbox{with }\ x_1<x_2.
\EQN
\BT\label{PreThm1}
Let $X_u(t),t\in \R$ be a family of centered Gaussian processes with  variance functions $\sigma^2_u$ and correlation functions $r_u$.
If  \textbf{A1}-\textbf{A3} are satisfied with $\Delta(u)=[x_1(u),x_2(u)]$, and $f\in C_0^*([x_1,x_2])$, then for $M_u$ satisfying ${M_u}\sim{u}, u\rw\IF $, we have
\BQN
\renewcommand\arraystretch{1.3}
\pk{\sup_{t\in\Delta(u)} X_u(t_u+t)>M_u}\sim C \LT(u^{\lambda}\orh(u^{-1})\RT)^{-\mathbb{I}_{\{\eta=\IF\}}} \Psi(M_u), \quad u\to \IF,
\EQN
 where
\BQN\label{C}
C=\left\{
\begin{array}{ll}
\mathcal{H}_{\alpha}\int_{x_1}^{x_2}e^{-f(t)}dt,&\hbox{if} \ \ \eta=\IF,\\
\mathcal{P}_{\alpha,\eta}^{f}[x_1,x_2],& \hbox{if} \ \ \eta\in(0,\IF),\\
\sup_{t\in[x_1,x_2]}e^{-f(t)},&  \hbox{if} \ \ \eta=0,
\end{array}
\right.
\EQN
and  $\mathcal{P}_{\alpha,\eta}^{f} (-\IF,\IF)\in(0,\IF)$.
\ET
\begin{remark}\label{re23}
\cJ{Let $\alpha\in (0,2], a>0$ be given. If $f\in C^*_0([x_1,x_2])$ for $ x_1,x_2,y\inr, x_1< x_2$, as shown in Appendix, we have, with $f_y(t):= f(y+t),t\inr$}
	\BQN\label{prop1} \mathcal{P}_{\alpha,a}^{f}[x_1,x_2]=\mathcal{P}_{\alpha,a}^{f_y}[x_1-y,x_2-y], \quad
	\mathcal{P}_{\alpha,a}^{f}[x_1,\IF)=\mathcal{P}_{\alpha,a}^{f_y}[x_1-y,\IF).
	\EQN
In particular, if  $f(t)=ct, c>0$, then  for any  $x\inr $
\begin{align*}
\mathcal{P}_{\alpha,a}^{ct}[x,\IF)
=\mathcal{P}_{\alpha,a}^{cx+ct}[0,\IF)
=e^{-cx}\mathcal{P}_{\alpha,a}^{ct}[0,\IF).
\end{align*}

\end{remark}
Next, for any fixed $T\in(0,\IF)$, in order to analyse $p_T(u)$  we shall suppose that:\\
\textbf{A1'}: For all large $u$, \cJ{$\sigma_u(t)$ attains its maximum  over $[0,T]$ at a unique} point $t_u$ such that
$$
\sigma_u(t_u)=1\quad\hbox{and}\quad \lim_{u\to\IF} t_u=t_0\in[0,T].
$$
\textbf{A4}: For all $u$ large enough
$$
\inf_{t\in[0,T]\setminus(t_u+\Delta(u))}\frac{1}{\sigma_u(t)}\geq1+\frac{p(\ln u)^q}{u^2}
$$
holds for some constants $p>0,q>1$.

\textbf{A5}: For some positive constants $G, \varsigma>0$
$$\E{( \overline{X}_u(t)-\overline{X}_u(s))^2}\leq G|t-s|^\varsigma$$
holds for all $s,t\in\{x\in[0,T]:\sigma(x)\neq 0\}$ and $\overline{X}_u(t)=\frac{X_u(t)}{\sigma_u(t)}$.\\
Below we define for $\lambda$ given in {\bf A2} and $\nu, d$ positve
\BQN\label{Du}
\quad \quad \quad \quad \quad \Delta(u)=
\LT\{
\begin{array}{ll}
	{[0,\delta_u]} & \  \text{if} \  \ t_u\equiv 0,\\
	{[-t_u,\delta_u]}, & \  \text{if} \  \ t_u\sim du^{-\nu}\ \text{and}\ \nu\geq\lambda,\\
	{[-\delta_u,\delta_u]}, & \  \text{if}\ \ t_u\sim du^{-\nu}\ \text{or}\ T-t_u\sim du^{-\nu}\ \text{when\ } \nu<\lambda,\ \text{or}\ t_0\in(0,T),\\
	{[-\delta_u,T-t_u]}, & \  \text{if} \  \ T-t_u\sim du^{-\nu}\ \text{and}\ \nu\geq\lambda,\\
	{[-\delta_u,0]} & \  \text{if} \  \ t_u= T,
\end{array}
\RT.
\EQN
where  $\delta_u=\LT(\frac{(\ln u)^q}{u}\RT)^\lambda$ with $q$ given in {\bf A4}.
\BT\label{PreThm2}
Let $X_u(t),t\in [0,T]$ be a family of centered Gaussian processes with   variance functions $\sigma^2_u$ and correlation functions $r_u$.
Assume that  \textbf{A1'},\textbf{A2}-\textbf{A5} are satisfied with $\Delta(u)=[c_1(u),c_2(u)]$ given in \eqref{Du} and
 $$ \lim_{u\rw\IF}c_i(u)u^{\lambda} = x_i\in [- \IF, \IF], i=1,2, \quad x_1 < x_2.$$
 If $f\in C^*_0([x_1,x_2])$, then for $M_u$ suc that $\limit{u} M_u/u=1$ we have
\BQN
\renewcommand\arraystretch{1.3}
\pk{\sup_{t\in[0,T]} X_u(t)>M_u}\sim C \LT(u^{\lambda}\orh(u^{-1})\RT)^{-\mathbb{I}_{\{\eta=\IF\}}} \Psi(M_u), \quad u\to \IF,
\EQN
where $C$ is the same as in \eqref{C} if $\eta\in(0,\IF]$ and $C=1$ if $\eta=0$.
\ET
\begin{remark}
% i)  When $\lambda=2/\beta$, $f(t)=b\abs{t}^\beta$ with $b,\beta>0$, and $\rho^2(t)=a\abs{t}^\alpha$ with $a>0,\alpha\in(0,2]$, the result of \netheo{PreThm1} reduces to the classical Piterbarg-Prisja{\v{z}}njuk
%	result, see \cite{PP78}.

  \cJ{\netheo{PreThm2}  generalises both  \cite{PP78}[Theorem 1] and  \cite{DHJ13a}[Theorem 4.1]. }

\end{remark}

\section{Applications}

%\cLa{Next, we focus our interest on the asymptotics of $\pk{\sup_{t\in[0,T]}(X(t)+g(t))>u}$, as $u\rw\IF$.}
\subsection{Locally stationary Gaussian processes with trend}
In this section we consider the asymptotics of %$$\pk{\sup_{t\in[0,T]}\LT(X(t)+g(t)\RT)>u}$$ as $u\rw\IF$, where
\eqref{Ab1} for $X(t),t\in[0,T]$ a centered  locally stationary Gaussian process with  unit  variance  and correlation function  $r$ satisfying
 \BQN \label{stationaryR0}
\lim_{h\to 0} \sup_{t\in [0,T]} \Abs{ \frac{ 1- r(t,t+h) }{a(t) \abs{h}^{\alpha}}} =1
% + o(\abs{h}^{\alpha}), \  \ h \to 0
\EQN
with  $\alpha \in (0,2]$,  $a(\cdot)$  a  positive continuous function  on $[0,T]$ and further
% Locally stationary Gaussian processe are introduced and studied in \cite{Ber74}.
\BQN \label{stationaryR2}
 r(s,t)< 1, \  \forall s,t\in[0,T]\ \mathrm{and\ } s \not=t.
\EQN

We refer to e.g., \cite{Ber74, Berman92, Hus90, Pit96,ChengD}  for results on  locally stationary  Gaussian processes. Extensions of this class to $\alpha(t)$-locally stationary processes are discussed in \cite{MR2462286,MR3413855, Long}.

Regarding the \cLa{continuous} trend function $g$, we define $g_m=\max_{t\in [0,T]} g(t)$  and set
$$H:=\LT\{s\in[0,T]:g(s)=g_m\RT\}.$$	
%	If $g_m$ is attained at some point $t_0 \in [0,T]$, then it is important to know the behaviour of $g_m -g(t)$ for $t$ close to $t_0$.
%For notational simplicity we shall assume that $t_0= argmax_{t\in [0,T]} g(t)$ and
%Without loss of generality, we shall assume
%$$ g_m=0.$$
	
 Set below, for any $t_0\in[0,T]$
\begin{align}\label{II}
Q_{t_0}=1+\mathbb{I}_{\{t_0\in(0,T)\}},\ \ \ \pa= \LT\{
\begin{array}{ll}
-\IF,& \hbox{if} \  t_0\in(0,T),\\
0,& \hbox{if} \  t_0=0\ \hbox{or}\ t_0=T.
\end{array}
\RT.
\end{align}

\begin{prop}\label{Thm3}
\kd{Suppose that \eqref{stationaryR0} and \eqref{stationaryR2} hold for
a centered locally stationary Gaussian process $X(t), t\in[0,T]$}
and let $g:[0,T]\to \R$ be a \cLa{continuous} function.  \\
{i)} If  $H=\{t_0\}$ and \eqref{eq:gtt0} holds, then as $u\rw\IF$
\BQN\label{ei}
\pk{\sup_{t\in[0,T]}(X(t)+g(t)) >u}\sim  C_{t_0} u^{(\frac{2}{\alpha}-\frac{1}{\gamma})_{+}}\Psi\LT( u-g_m\RT),
 \EQN
where (set with  $a=a(t_0)$)
$$ C_{t_0}= \left\{
\begin{array}{ll}
\qa a^{1/{\alpha}}c^{-1/\gamma}\Gamma(1/\gamma+1)\mathcal{H}_{\alpha},&\hbox{if} \ \ \alpha<2\gamma,\\
\mathcal{P}_{\alpha,a}^{c|t|^\gamma}[\pa,\IF),&  \hbox{if} \ \ \alpha=2\gamma,\\
1,&  \hbox{if} \ \ \alpha>2\gamma.
\end{array}
\right.
$$
{ii)} If $H=[A,B]\subset[0,T]$ with $0\leq A<B\leq T$, then as $u\to \IF$
\BQNY
\pk{\sup_{t\in[0,T]}(X(t)+g(t)) >u}\sim \mathcal{H}_{\alpha}\int_{A}^B(a(t))^{1/\alpha}dt u^{\frac{2}{\alpha}}\Psi\LT(u-g_m\RT).
\EQNY
 \end{prop}
\begin{remarks} \label{glass}
i) If $H=\kd{\{t_1 \ldot t_n\}}$, then as mentioned in \cite{Pit96}, the tail \kd{distribution}
of the \kd{corresponding} supremum is easily obtained assuming that for each $t_i$ the assumptions of \neprop{Thm3} statement {i)} hold, implying that
$$\pk{\sup_{t\in[0,T]}(X(t)+g(t)) >u} \sim \Bigl(\sum_{j=1}^n C_{t_j} \Bigr)  u^{(\frac{2}{\alpha}-\frac{1}{\gamma})_{+}}\Psi\LT( u-g_m\RT), \quad u\to \IF.
$$

ii) The \kd{novelty} of \neprop{Thm3} statement {i)} is that  for the trend function $g$ only a polynomial local behavior around $t_0$ is assumed. In the literature so far only the case that \eqref{eq:gtt0} holds with $\gamma= 2$ has been considered  (see \cite{PitSta01}).
% Hence, for  $\gamma \le 1$, since $\alpha \in (0,2]$, both the case $\alpha=2 \gamma$ in \neprop{Thm3}  and the case $\alpha > 2 \gamma$  are completely new.

%\cJ{ iii) For general case where $g_m\neq 0$, the results in Proposition \ref{Thm3} still hold, but with $u$ replaced by $u-g_m.$}

\cLa{ iii) \kd{By the proof of Proposition \ref{Thm3} statement i),}
if $g(t)$ \kd{is} a measurable function which is continuous in a neighborhood of $t_0$
and smaller \kd{than} $g_m-\vn$ for some $\vn>0$ in the rest part over $[0,T]$, then the results still hold.}
\end{remarks}

We present below the approximation of the conditional passage time $\tau_u|\tau_u\leq T$ with $\tau_u$ defined 	in \eqref{tu}.
\begin{prop}\label{ruintimelocal}
\kd{Suppose that \eqref{stationaryR0} and \eqref{stationaryR2} hold for
a centered locally stationary Gaussian process $X(t), t\in[0,T]$.
Let $g:[0,T]\to \R$ be a continuous function,
$H=\{t_0\}$ and \eqref{eq:gtt0} holds.
\\
%Under the assumptions and notation in of \neprop{Thm3} statement {i)}, when
i) If $t_0\in[0,T)$, then for any  $x\in (\pa,\IF)$}
\begin{align*}
\renewcommand\arraystretch{1.7}
\pk{u^{1/\gamma}(\tau_u-t_0)\leq x\big| \tau_u\leq T}\sim
\LT\{
\begin{array}{ll}
\frac{{\gamma}c^{1/\gamma}\int_{\pa}^{x} e^{-c\abs{t}^{\gamma}}dt}{ {\qa}\Gamma(1/\gamma)},& \text{if}\ \alpha<2\gamma,\\
\frac{ \mathcal{P}_{\alpha,a}^{c\abs{t}^{\gamma}}[\pa,x]}{ \mathcal{P}_{\alpha,a}^{c\abs{t}^{\gamma}}[\pa,\IF)}, &  \text{if}\ \alpha=2\gamma,\\
\sup_{t\in[\pa,x]}e^{-c\abs{t}^{\gamma}},&  \text{if}\ \alpha>2\gamma,
\end{array}
\RT.
\end{align*}
\kd{ii) If $t_0=T$, then for any  $x\in (-\IF,0)$}
\begin{align*}
\renewcommand\arraystretch{1.7}
\pk{u^{1/\gamma}(\tau_u-t_0)\leq x\big| \tau_u\leq T}\sim
\LT\{
\begin{array}{ll}
\frac{{\gamma}c^{1/\gamma}\int_{-x}^{\IF} e^{-c\abs{t}^{\gamma}}dt}{\Gamma(1/\gamma)},& \text{if}\ \alpha<2\gamma,\\
\frac{ \mathcal{P}_{\alpha,a}^{c\abs{t}^{\gamma}}[-x,\IF)}{ \mathcal{P}_{\alpha,a}^{c\abs{t}^{\gamma}}[0,\IF)}, &  \text{if}\ \alpha=2\gamma,\\
 e^{-c\abs{x}^{\gamma}},&  \text{if}\ \alpha>2\gamma.
\end{array}
\RT.
\end{align*}
\end{prop}
\COM{\begin{remark}
In the above corollary we used the fact that
\BQNY
\mathcal{P}_{\alpha,a}^{c\abs{t}^{\gamma}}[-x,\IF)=\mathcal{P}_{\alpha,a}^{c\abs{t}^{\gamma}}(-\IF, x], \ \ x\in(-\IF,0).
\EQNY
\end{remark}}

%The next example extends \cite{GPPit2015}[Lemma 2].
%Note that the Ornstein-Uhlenbeck process with trend is a special case with $\alpha=1$.
\begin{example}\label{exam2}
	Let $X(t), t\in[0,T]$ be a centered stationary Gaussian process with  unit variance and  correlation function $r$ that satisfies $r(t)=1-a |t|^\alpha(1+o(1)),\ t\rw0$ for some $a>0$, $\alpha\in(0,2]$, and $r(t)<1$, for all $t\in(0,T]$. \cJ{Let $\tau_u$ be defined as in \eqref{tu} with $g(t)=-ct, c>0$}. Then we have
	\BQNY
	\pk{\max_{t\in[0,T]}\LT(X(t)-ct\RT)>u }\sim u^{(\frac{2}{\alpha}-1)_{+}}\Psi(u )
	\LT\{\begin{array}{ll}
		c^{-1}a^{1/\alpha}\mathcal{H}_\alpha,& \ \alpha\in(0,2),\\
		\mathcal{P}_{\alpha,a}^{ct}[0,\IF),& \ \alpha=2,
	\end{array}
	\RT.
	\EQNY
and for any $x$ positive
\BQNY
\pk{u\tau_u\leq x \Big| \tau_u\leq T}\sim
	\LT\{\begin{array}{ll}
		1-e^{-cx},& \ \alpha\in(0,2),\\
		\frac{\mathcal{P}_{\alpha,a}^{ct}[0,x]}{\mathcal{P}_{\alpha,a}^{ct}[0,\IF)},& \ \alpha=2.
	\end{array}
	\RT.
\EQNY
\COM{Similarly,  we have
	\BQNY
	\pk{\max_{t\in[0,T]}\LT(X(t)-ct^2\RT)>u }\sim \frac{\sqrt{\pi} }{2} c^{-\frac{1}{2}}a^{1/\alpha}\mathcal{H}_\alpha u^{\frac{2}{\alpha}-\frac{1}{2}}\Psi(u ),
	\EQNY
and
\BQN
\pk{u^{1/2}\tau_u\leq x \Big| \tau_u\leq T}\sim
 %2\sqrt{\frac{c}{\pi}}\int_{0}^xe^{-ct^2}dt=
 2\Phi(\sqrt{2c}x)-1,	
\ \ x\in(0,\IF).
\EQN}

\end{example}

\begin{example}\label{exam3}
	Let $X(t), t> 0$ be a standardized fBm, i.e., $X(t)= B_{\alpha}(t) / t^{\alpha/2}$ \cJ{with $B_{\alpha}$ \kd{an} fBm.}  %with self-similarity index $\alpha/2 \in(0,1]$.
	Let  $c, T$ be positive constants.
Then for any $n\in \N$, we have
	\BQNY
	\pk{\max_{t\in[T,(n+1)T]}\LT(X(t)+c \sin\LT(\frac{2\pi t}{T}\RT)\RT)>u }
\sim\LT(\sum_{j=1}^n a_j^{\frac{1}{\alpha}}\RT)\mathcal{H}_{\alpha}
	\frac{T}{\sqrt{2c\pi}} u^{\frac{2}{\alpha}-\frac{1}{2}}\Psi(u-c),
	\EQNY
	where $a_j=\frac{1}{2}\LT(\frac{(4j+1)T}{4}\RT)^{-\alpha},j=1,\ldots,n$.
\end{example}
\subsection{Non-stationary Gaussian processes with trend}
In this section we consider the asymptotics of %$$\pk{\sup_{t\in[0,T]}\LT(X(t)+g(t)\RT)>u}$$ as $u\rw\IF$, where
\eqref{Ab1} for  $X(t),t\in[0,T]$ a centered Gaussian process with non-constant \kd{variance function $\sigma^2$}. % and  correlation function $r$.
Define below %$X_g(t):=X(t)+g(t)$
whenever  $\sigma(t)\neq 0$
$$\overline{X}(t):=\frac{X(t)}{\sigma(t)}, \quad t\in [0,T],$$
and set for  \kd{a continuous function $g$}
\BQN\label{eq:Mtu}
m_u(t):=\frac{\sigma(t)}{1-g(t)/u},\ \ \ t\in[0,T], \quad u>0.
\EQN

\begin{prop}\label{MainThm1}
\cJ{Let  $X$ and $g$ be as above. Assume that $t_u=argmax_{t\in [0,T]}m_u(t)$ is unique with $\lim_{u\to\IF}t_u= t_0$ and $\sigma(t_0)=1$.}   Further, we suppose that  \textbf{A2}-\textbf{A5} are satisfied with $\sigma_u(t)=\frac{m_u(t)}{m_u(t_u)}$, $r_u(s,t)=r(s,t)$, $\overline{X}_u(t)=\overline{X}(t)$ and $\Delta(u)=[c_1(u),c_2(u)]$ given in (\ref{Du}).
 If  in \textbf{A2} $f \in C^*_0([x_1,x_2])$   and
  $$ \lim_{u\rw\IF}c_i(u)u^{\lambda} = x_i\in [- \IF, \IF], i=1,2, \quad x_1 < x_2,$$
 then we have
\BQN\label{ei2}
%p_T(u)
\pk{\sup_{t\in[0,T]}\Xgt >u}\sim  C \LT(u^\lambda\orh(u^{-1})\RT)^{-\mathbb{I}_{\{\eta=\IF\}}}\Psi\LT(\frac{u-g(t_u)}{\sigma(t_u)}\RT), \quad u\to \IF,
\EQN
where $C$ is the same as in \eqref{C} when $\eta\in(0,\IF]$ and $C=1$ when $\eta=0$.
\end{prop}
\begin{remarks}
i) \neprop{MainThm1} extends  \cite{PP78}[Theorem 3]   and the results of \cite{DebRol02} where \eqref{Ab1} was analyzed for special  $X$ with stationary increments and special trend function $g$.

\cJ{
ii) The assumption that $\sigma(t_0)= 1$ is not essential in the proof. In fact, for the general case where $\sigma(t_0)\neq 1$ we have that \eqref{ei2} holds with
\BQNY
C=\left\{
\begin{array}{ll}
\sigma_0^{-\frac{2}{\alpha}}\mathcal{H}_{\alpha}\int_{x_1}^{x_2}e^{-\sigma_0^{-2}f(t)}dt,&\hbox{if} \ \ \eta=\IF,\\
\mathcal{P}_{\alpha,\sigma_0^{-2}\eta}^{\sigma_0^{-2}f}[x_1,x_2],& \hbox{if} \ \ \eta\in(0,\IF),\\
1,&  \hbox{if} \ \ \eta=0,
\end{array}
\right. \ \ \ \ \sigma_0=\sigma(t_0).
\EQNY}
\end{remarks}
\begin{prop} \label{kA}
Under the notation and assumptions of \neprop{MainThm1} \cJ{without assuming {\bf A3,A5},} if $X$ is  differentiable in the mean square \cJ{sense} such that
$$r(s,t)<1, s\neq t, \quad \E{X'^2(t_0)}>\cJ{\sigma'^2(t_0)},$$
and $\E{X'^2(t)}-\sigma'^2(t)$ is continuous in a neighborhood of $t_0$,
  then \cJ{\eqref{ei2}} holds with
  $$\alpha=2,\ \ \  \rho^2(t)=\frac{1}{2}\LT(\E{{X'}^2(t_0)}-\sigma'^2(t_0)\RT)t^2.$$

\end{prop}

\COM{\BK\label{cor1}
 For $X(t)$ and $g(t)$ above, suppose that $X(t)$ is differentiable in mean square and \textbf{A1'},\textbf{A2},\textbf{A4} and \textbf{A5} are satisfied with $\sigma_u(t)=\frac{m_u(t)}{m_u(t_u)}$ and $\Delta(u)$ in (\ref{Du}) with $[x_1,x_2]=\lim_{u\rw\IF}u^{\lambda}\Delta(u)$. Moreover, we assume that $r(s,t)<1, s\neq t$ and $\E{X'^2(t_0)}-\sigma'^2(t_0)>0$.
Set $\sigma_0=\sigma(t_0)$, then we have
\BQN
\pk{\sup_{t\in[0,T]}\Xgt >u}\sim  \Psi\LT(\frac{u-g(t_u)}{\sigma(t_u)}\RT)\left\{
\begin{array}{ll}
\frac{u^{1-\lambda}\eta^{1/2}}{\sqrt{\pi}\sigma_0}\int_{x_1}^{x_2} e^{-w(t)}dt,&\hbox{if} \ \ \lambda<1,\\
 \mathcal{P}_{2,\sigma_0^{-2}\eta}^{w}[x_1,x_2],& \hbox{if} \ \  \lambda=1,\\
1,&  \hbox{if} \ \ \lambda>1,
\end{array}
 \right.
\EQN
where $w(t)$ is the same as in \netheo{MainThm1} \cl{and} $\eta=\frac{\E{{X'}^2(t_0)}-\sigma'^2(t_0)}{2\sigma_0^2}$. % and $\alpha=2$.
\EK}

The next result is an extension of a classical theorem concerning the extremes of non-stationary Gaussian processes \cJ{discussed in the Introduction,} see \cite{Pit96}[Theorem D.3].
\begin{prop}\label{Thm2}
\EHc{Let $X(t),t\in [0,T]$ be a centered Gaussian process with correlation function $r$ and variance function
$\sigma^2$ such that $t_0=argmax_{t\in [0,T]}\sigma(t)$ is unique with $\sigma(t_0)=\sigma>0$. Suppose  that  $g$ is a bounded measurable  function being continuous in a neighborhood of $t_0$ such that \eqref{eq:gtt0} holds. If further \eqref{assump-corre} is satisfied}, then
\BQN\pk{\sup_{t\in[0,T]}\Xgt >u}\sim C_{0} u^{(\frac{2}{\alpha}-\frac{2}{\beta^*})_{+}} \Psi\LT(\frac{u-g(t_0)}{\sigma}\RT), %\quad \sigma:=\sigma(t_0),
\EQN
where
$\beta^*=\min(\beta,2\gamma)$,
$$C_{0}=\LT\{
\begin{array}{ll}
\sigma^{-2/\alpha} a^{1/\alpha}\mathcal{H}_\alpha\int_{\pa}^\IF e^{-f(t)}dt,& \text{if}\ \alpha<\beta^*,\\
\mathcal{P}_{\alpha,\sigma^{-2} a}^{f}[\pa,\IF), &  \text{if}\ \alpha=\beta^*,\\
1,&  \text{if}\ \alpha>\beta^*,
\end{array}
\RT.
$$
with
$f(t)=\frac{b}{\sigma^3}|t|^{\beta} \mathbb{I}_{\{\beta=\beta^*\}}+\frac{c}{\sigma^2}|t|^\gamma \mathbb{I}_{\{2\gamma=\beta^*\}}$
and $\pa$  \kd{defined} in \eqref{II}.
\end{prop}

\begin{prop}\label{ruintime}
i) Under the conditions and notation of \neprop{MainThm1}, for any  $x\in [x_1,x_2]$  we have
\BQN
\renewcommand\arraystretch{1.7}
\pk{u^{\lambda}(\tau_u-t_u)\leq x\big| \tau_u\leq T}\sim  \left\{
\begin{array}{ll}
\frac{\int_{x_1}^{x} e^{-f(t)}dt}{\int_{x_1}^{x_2} e^{-f(t)}dt},&\hbox{if} \ \ \eta=\IF,\\
\frac{ \mathcal{P}_{\alpha,\eta}^{f}[x_1,x]}{ \mathcal{P}_{\alpha,\eta}^{f}[x_1,x_2]},& \hbox{if} \ \  \eta\in(0,\IF),\\
\sup_{t\in[x_1,x]}e^{-f(t)},&  \hbox{if} \ \ \eta=0.
\end{array}
 \right.
\EQN
ii) Under the conditions and notation of \neprop{Thm2},  \kd{if} $t_0\in[0,T)$, then for  $x\in (\pa,\IF)$
\begin{align*}
\renewcommand\arraystretch{1.7}
\pk{u^{2/{\beta^*}}(\tau_u-t_0)\leq x\big| \tau_u\leq T}\sim
\LT\{
\begin{array}{ll}
\frac{\int_{\pa}^{x} e^{-f(t)}dt}{\int_{\pa}^{\IF} e^{-f(t)}dt},& \text{if}\ \alpha<\beta^*,\\
\frac{ \mathcal{P}_{\alpha,a}^{f}[\pa,x]}{ \mathcal{P}_{\alpha,a}^{f}[\pa,\IF)}, &  \text{if}\ \alpha=\beta^*,\\
\sup_{t\in[\pa,x]}e^{-f(t)},&  \text{if}\ \alpha>\beta^*,
\end{array}
\RT.
\end{align*}
and if  $t_0=T$, then for   $x\in (-\IF,0)$
\begin{align*}
\renewcommand\arraystretch{1.7}
\pk{u^{2/\beta^*}(\tau_u-t_0)\leq x\big| \tau_u\leq T}\sim
\LT\{
\begin{array}{ll}
\frac{\int_{-x}^{\IF} e^{-f(t)}dt}{\int_{0}^{\IF} e^{-f(t)}dt},& \text{if}\ \alpha<\beta^*,\\
\frac{ \mathcal{P}_{\alpha,a}^{f}[-x,\IF)}{ \mathcal{P}_{\alpha,a}^{f}[0,\IF)}, &  \text{if}\ \alpha=\beta^*,\\
e^{-f(x)},&  \text{if}\ \alpha>\beta^*.
\end{array}
\RT.
\end{align*}
\end{prop}

\begin{example}\label{exam1}
	Let $X(t)=B(t)-tB(1), t\in[0,1]$, where $B(t)$ is a standard Brownian motion \kd{and
suppose that $\tau_u$ is defined by \eqref{tu} with $g(t)=-ct$. Then}
	\BQN\label{ex1eq1}
	\pk{\sup_{t\in[0,1]}(X(t)-ct)>u }\sim e^{-2(u^2+cu)},
	\EQN
    \BQNY
    \pk{u\LT(\tau_u-\frac{u}{c+2u}\RT)\leq x \Big| \tau_u\leq 1}\sim\Phi(4x),\ \ x\in(-\IF,\IF).
    \EQNY
\kd{We} note that 	according to \cite{AsymBB2003}[Lemma 2.7], the result in \eqref{ex1eq1}
is actually exact, i.e. for any $u>0$,\\
$\pk{\sup_{t\in[0,1]}(X(t)-ct)>u }= e^{-2(u^2+cu)}.$

\kd{Now, let $T=1/2$. It appears that the asymptotics in this case is} different, i.e.,
	\BQN\label{ex1eq2}
	\pk{\sup_{t\in[0,1/2]}(X(t)-ct)>u }\sim \Phi(c) e^{-2(u^2+cu)},
	\EQN
and
    \BQNY
    \pk{u\LT(\tau_u-\frac{u}{c+2u}\RT)\leq x \Big| \tau_u\leq \frac{1}{2}}\sim
    \frac{\Phi(4x)}{\Phi(c)},\ \ x\in(-\IF,c/4].
    \EQNY
	Similarly, we have
	\BQN\label{ex1eq3}
	\pk{\sup_{t\in[0,1]}\LT(X(t)+\frac{c}{2}-c\LT|t-\frac{1}{2}\RT|\RT)>u }\sim 2\Psi(c)e^{-2(u^2-cu)}
	\EQN
     and
    \BQNY
    \pk{u\LT(\tau_u-\frac{1}{2}\RT)\leq x \Big| \tau_u\leq 1}\sim \frac{\int_{-\IF}^{4x}e^{-\frac{(|t|+c)^2}{2}}d t}{2\sqrt{2\pi}\Psi(c)},\ \ x\in(-\IF,\IF).
    \EQNY
\end{example}

\cJ{We conclude this section \kd{with}} an application of \neprop{MainThm1} to the calculation
of the ruin probability of a Brownian motion risk model \cJ{with constant force of interest} over infinite-time horizon.
\subsection{Ruin probability in Gaussian risk model}\label{ruin}
\kd{Consider risk reserve process $U(t)$, with interest rate $\delta$, modeled by} %$ U(t),t\ge0$ given as
$$U(t)=ue^{\delta t}+c \int_0^t e^{\delta (t-v)}dv-\sigma \int_0^t e^{\delta(t-v)}dB(v), \quad t\ge0,$$
where $c,\delta, \sigma$ are some positive constants and $B$ is a standard Brownian motion.
The corresponding ruin probability over infinite-time horizon is defined as
$$ p(u)= \pk{\inf_{t\in[0,\IF)}U(t)<0}. $$
For this model we also define the ruin time   $\tau_u=\inf\{t\geq0:U(t) < 0\}$. Set below
$$h(t)=\frac{\delta}{\sigma^2}\LT(\sqrt{t+r^2}-r\RT)^2,\quad  t\in [0,\IF),\quad r=\frac{c}{\delta}.$$

We present next approximations of the ruin probability and the conditional ruin time $\tau_u|\tau_u<\IF$ as $u\to \IF$.
\BS\label{exmodel} As $u\to \IF$
\BQN\label{ruinex}
p(u)\sim
\mathcal{P}_{1,\delta/\sigma^2}^{h}\LT[-r^2,\IF\RT)\Psi\LT(\frac{1}{\sigma}\sqrt{2\delta u^2+4cu}\RT)
\EQN
and for $x\in(-r^2,\IF)$
\BQNY
\pk{u^{2}\LT(e^{-2\delta\tau_u}-\LT(\frac{c}{\delta u+c}\RT)^2\RT)\leq x\big| \tau_u<\IF}\sim \frac{\mathcal{P}_{1,\delta/\sigma^2}^{h}\LT[-r^2,x\RT]}{\mathcal{P}_{1,\delta/\sigma^2}^{h}\LT[-r^2,\IF\RT)}.
\EQNY
\ES
%\end{example}
\begin{remark}\label{pep}
According to \cite{RPwCA1977} (see also \cite{ADARP1975})  we have
\BQN\label{CC}
\pk{\inf_{t\in[0,\IF)}U(t)<0}=\Psi\LT(\frac{\sqrt{2\delta}}{\sigma}\LT(u+r\RT)\RT)\Big/\Psi\LT(\frac{\sqrt{2}c}{\sigma\sqrt{\delta}}\RT).
\EQN
By \eqref{ruinex} and \eqref{prop1}
\BQNY
\pk{\inf_{t\in[0,\IF]} U(t)<0}
&\sim&\E{\sup_{t\in[-r^2,\IF)}
\exp\LT(\sqrt{\frac{2\delta}{\sigma^2}}B(t)-\frac{\delta}{\sigma^2}\LT(\sqrt{t+r^2}-r\RT)^2
-\frac{\delta}{\sigma^2}|t|\RT)}\Psi\LT(\frac{1}{\sigma}\sqrt{2\delta u^2+4cu}\RT)\\
&\sim&\E{\sup_{t\in[-\frac{c^2}{\sigma^2\delta},\IF)}
\exp\LT(\sqrt{2}B(t)-\LT(t+\frac{c^2}{\sigma^2\delta}\RT)+\frac{2c}{\sigma\sqrt{\delta}}\sqrt{t+\frac{c^2}{\sigma^2\delta}}-|t|\RT)}
\Psi\LT(\frac{\sqrt{2\delta}}{\sigma}\LT(u+r\RT)\RT)\\
%&\sim&\E{\sup_{t\in[0,\IF)}
%\exp\LT(\sqrt{2}B\LT(t-\frac{c^2}{\sigma^2\delta}\RT)-\LT(t-\frac{c}{\sigma\sqrt{\delta}}\RT)^2-\LT|t-\frac{c^2}{\sigma^2\delta}\RT|\RT)}
%\Psi(M_u)\\
&=&\E{\sup_{t\in[0,\IF)}
\exp\LT(\sqrt{2}B\LT(t\RT)-2t+\frac{2c}{\sigma\sqrt{\delta}}\sqrt{t}\RT)}
\Psi\LT(\frac{\sqrt{2\delta}}{\sigma}\LT(u+r\RT)\RT),
\EQNY
which combined with \eqref{CC} implies that
\BQN
\E{\sup_{t\in[0,\IF)}
\exp\LT(\sqrt{2}B\LT(t\RT)-2t+\frac{2c}{\sigma\sqrt{\delta}}\sqrt{t}\RT)}
= \LT(\Psi\LT(\frac{\sqrt{2}c}{\sigma\sqrt{\delta}}\RT)\RT)^{-1}.
\EQN
\end{remark}

\section{Proofs}
\kd{In the proofs presented in this section} $\mathbb{C}_i, i\in\N$ are some positive constants which \cJ{may be different from line to line}.\\
 We first give two preliminary lemmas, which play an important role in the proof of \netheo{PreThm1}.
\BEL\label{pan2} Let $\xi(t), t\in\R$ be a centered stationary Gaussian process with  unit variance and correlation function $r$ satisfying
\BQN\label{rp}
1-r(t)\sim a\rho^2(\abs{t}),\ \ \ t\to 0,
\EQN
with $a>0$, and $\rho\in\mathcal{R}_{\alpha/2}$, $\alpha\in(0,2]$. Let $f$ be a continuous function, $K_u$ be a family of index sets and
\BQNY
Z_{u}(t):=\frac{\xi(\orh(u^{-1})t)}{1+u^{-2}f(\orh(u^{-1})u^{\lambda}t)},\quad  t\in[S_1,S_2],
\EQNY
where $\lambda>0$ and $-\IF<S_1<S_2<\IF$.
If $M_k(u), k\in K_u$ is such that
\BQN\label{Mku}
\lim_{u\to\IF}\sup_{k\in K_u}\LT|\frac{M_k(u)}{u}-1\RT|=0,
\EQN
then we have
\BQN\label{qq}
\lim_{u\rw\IF}\sup_{k\in K_u}\LT|\frac{1}{\Psi(M_k(u))}\pk{\sup_{t\in [S_1,S_2]}Z_{u}(t)>M_k(u)}-\mathcal{R}_\eta^f[S_1,S_2]\RT|=0,
\EQN
where
\begin{align*}
\mathcal{R}_\eta^f[S_1,S_2]:=
\E{\sup_{t\in[S_1,S_2]}e^{\sqrt{2a}B_\alpha(t)-a|t|^\alpha-f(\eta^{-1/\alpha}t)}}
=\LT\{
\begin{array}{ll}
\mathcal{H}_{\alpha}[a^{1/\alpha}S_1,a^{1/\alpha}S_2]&\  \ \cLa{f(\cdot)\equiv0,}\\
\mathcal{P}_{\alpha,a}^h[S_1,S_2]&\ \cLa{\hbox{otherwise}},
\end{array}
\RT.
\end{align*}
with $\eta:= \lim_{t\downarrow 0}\frac{\rho^2(t)}{t^{2/\lambda}}\in(0,\IF]$ and
\cLa{$h(t)= f(\eta^{-1/\alpha}t)$ for $\eta\in(0,\IF)$, $h(t)=f(0)$ for $\eta=\IF$.}
\EEL

\prooflem{pan2} We set $\eta^{-1/\alpha}=0$ if $\eta=\IF$.
The proof follows by checking the conditions of \cite{Uniform2016}[Theorem 2.1] \cLa{where the results still holds if we omit the requirements $f(0)=0$ and $[S_1,S_2]\ni 0$}.
By (\ref{Mku})
\BQNY
\lim_{u\rw\IF}\inf_{k\in K_u} M_k(u)=\IF.
\EQNY
\kd{By continuity of $f $ we have}
\BQN
\lim_{u\rw\IF}\sup_{k\in K_u,t\in[S_1,S_2]}\LT|M^2_k(u)u^{-2}f(\orh(u^{-1})u^{\lambda}t)-f(\eta^{-1/\alpha}t)\RT|=0.
\EQN
\kd{Moreover, \eqref{rp} implies}
\BQNY
\Var(\xi(\orh(u^{-1})t)-\xi(\orh(u^{-1})t'))=2-2r\LT(\abs{\orh(u^{-1})(t-t')}\RT)
\sim 2a\rho^2 \LT(\abs{\orh(u^{-1})(t-t')}\RT),\ u\rw\IF,
\EQNY
holds for $t,t'\in[S_1,S_2]$.
Thus
\BQN
\lim_{u\rw\IF}\sup_{k\in K_u}\sup_{t\neq t'\in[S_1,S_2]}\LT|M^2_k(u)\frac{\Var(\xi(\orh(u^{-1})t)-\xi(\orh(u^{-1})t'))}{ 2au^{2}\rho^2 \LT(\abs{\orh(u^{-1})(t-t')}\RT)}-1\RT|=0.
\EQN
Since $\rho^2\in\mathcal{R}_\alpha$ which satisfies the uniform convergence theorem (UCT) for regularly varying function, see, e.g., \cite{bingham1989regular}, i.e.,
\BQN\label{UCT}
\lim_{u\rw\IF}\sup_{t,t'\in[S_1,S_2]}\LT|u^{2}\rho^2 \LT(\abs{\orh(u^{-1})(t-t')}\RT)-\abs{t-t'}^\alpha\RT|=0,
\EQN
and further by the Potter's bound for $\rho^2$, see \cite{bingham1989regular} we have
\BQN\label{pott}
\limsup_{u\rw\IF}\underset{t\neq t'}{\sup_{t,t'\in[S_1,S_2]}}\frac{u^{2}\rho^2 \LT(\abs{\orh(u^{-1})(t-t')}\RT)}{\abs{t-t'}^{\alpha-\vn_1}}\leq \mathbb{C}_1\max\LT(\abs{S_1-S_2}^{\alpha-\vn_1},\abs{S_1-S_2}^{\alpha+\vn_1}\RT)<\IF,
\EQN
where $\vn_1\in(0,\min(1,\alpha))$.
We know that for $\alpha\in(0,2]$
\BQN\label{boundt}
\abs{\abs{t}^\alpha-\abs{t'}^\alpha}\leq\mathbb{C}_2\abs{t-t'}^{\alpha\wedge 1},\ t,t'\in[S_1,S_2].
\EQN
By \eqref{rp} for any small $\epsilon>0$, when $u$ large enough
\BQN
r(\orh(u^{-1})t)\leq1-\rho^2(\orh(u^{-1})\abs{t})(1-\epsilon),\quad
r(\orh(u^{-1})t)\geq1-\rho^2(\orh(u^{-1})\abs{t})(1+\epsilon)
\EQN
hold for $t\in[S_1,S_2]$,
then
by \eqref{Mku} for $u$ large enough
\BQN
&&\sup_{k\in K_u}\sup_{|t-t'|<\vn,t,t'\in [S_1,S_2]}M^2_k(u)\E{[\xi(\orh(u^{-1})t)-\xi(\orh(u^{-1})t')]\xi(0)}\nonumber\\
&& \leq\mathbb{C}_3 u^2\sup_{|t-t'|<\vn,t,t'\in [S_1,S_2]} \abs{r(\orh(u^{-1})t)-r(\orh(u^{-1})t')}\nonumber\\
&& \leq\mathbb{C}_3\sup_{|t-t'|<\vn,t,t'\in [S_1,S_2]} \LT(\abs{u^2\rho^2(\orh(u^{-1})\abs{t})-u^2\rho^2(\orh(u^{-1})\abs{t'})} + \epsilon\abs{u^2\rho^2(\orh(u^{-1})\abs{t})} +\epsilon\abs{u^2\rho^2(\orh(u^{-1})\abs{t'})}\RT) \nonumber\\
&&\leq \mathbb{C}_3 \sup_{|t-t'|<\vn,t,t'\in [S_1,S_2]} \LT(\abs{u^2\rho^2 \LT(\abs{\orh(u^{-1})(t)}\RT)-\abs{t}^\alpha}+\abs{u^2\rho^2 \LT(\abs{\orh(u^{-1})(t')}\RT)-\abs{t'}^\alpha}+\abs{\abs{t}^\alpha-\abs{t'}^\alpha}\RT.\nonumber\\
&&\quad\LT.+\mathbb{C}_4\epsilon\LT(\abs{t}^{\alpha-\vn_1}
+\abs{t'}^{\alpha-\vn_1}\RT)\RT)\label{expl1}\\
&&\leq \mathbb{C}_5 \vn^{\alpha\wedge 1}+\mathbb{C}_6\epsilon, \ u\rw\IF\label{expl2}\\
&&\rw 0,\vn\rw 0,\epsilon\rw0,\nonumber
\EQN
where in \eqref{expl1} we use \eqref{pott} and \eqref{expl2} follows from \eqref{UCT} and \eqref{boundt}.\\
Hence the proof follows from \cite{Uniform2016}[Theorem 2.1].
\QED

\BEL\label{pan3} Let $Z_u(s,t), (s,t)\in\R^2$ be a centered stationary Gaussian field with unit variance and correlation function $r_{Z_u}(\cdot,\cdot)$ satisfying
\BQN\label{rp1}
1-r_{Z_u}(s,t)=a u^{-2} \left(\Bigl\lvert\frac{s}{\orh(u^{-1})}\Bigl\lvert^{\alpha/2}+\Bigl\lvert\frac{t}{\orh(u^{-1})}\Bigl\lvert^{\alpha/2}\right),\ \ \ (s,t)\in\R^2,
\EQN
with $a>0$, $\rho^2\in\mathcal{R}_\alpha$ and $\alpha\in(0,2]$. Let $K_u$ be some index sets. Then, \cJ{for  $M_k(u), k\in K_u$ satisfying} \eqref{Mku}
and for any $S_1,S_2, T_1, T_2\geq0$ such that $\max(S_1,S_2)>0, \max(T_1,T_2)>0$, we have
\BQNY
\lim_{u\rw\IF}\sup_{k\in K_u}\LT| \frac{1}{\Psi( M_k(u))}\pk{\sup_{(s,t)\in D(u)}Z_u(s,t)>M_k(u)}-\mathcal{F}( S_1,S_2,T_1,T_2)\RT|=0,
\EQNY
where
 $D(u)=[-\orh(u^{-1})S_1,\orh(u^{-1})S_2]\times[-\orh(u^{-1})T_1,\orh(u^{-1})T_2]$ and
\BQNY
\mathcal{F}( S_1,S_2,T_1,T_2) %&:=& \E{\sup_{(s,t)\in [-S_1,S_2]\times[-T_1,T_2] }e^{c_0\sqrt{2a}(B_{\frac{\alpha}{2}}(s)+B'_{\frac{\alpha}{2}}(t))-ac_0^2(|s|^{\frac{\alpha}{2}}+|t|^{\frac{\alpha}{2}})}}\\
 = \mathcal{H}_{\alpha/2}[-a^{2/\alpha}S_1,a^{2/\alpha}S_2]
\mathcal{H}_{\alpha/2}[-a^{2/\alpha}T_1,a^{2/\alpha}T_2].
\EQNY
\EEL
\prooflem{pan3}
 The proof follows by checking the conditions  of \cite{KEP2015}[Lemma 5.3].\\
For $D=[-S_1,S_2]\times[-T_1,T_2]$ we have
\BQNY
\pk{\sup_{(s,t)\in D_u}Z_u(s,t)>M_k(u)}=
\pk{\sup_{(s,t)\in D}Z_u(\orh(u^{-1})s,\orh(u^{-1})t)>M_k(u)}.
\EQNY
Since by (\ref{rp1})
\BQNY
\Var(Z_u(\orh(u^{-1})s,\orh(u^{-1})t)-Z_u(\orh(u^{-1})s',\orh(u^{-1})t'))
&=&2-2r_{Z_u}\LT(\orh(u^{-1})(s-s'),\orh(u^{-1})(t-t')\RT)\\
&=&a u^{-2} \left(\abs{ s-s'}^{\alpha/2}+\abs{ t-t'}^{\alpha/2}\right)
\EQNY
we obtain
\BQN
\lim_{u\rw\IF}\sup_{k\in K_u}\sup_{(s,t)\neq (s',t')\in D}\LT|M^2_k(u)\frac{\Var(Z_u(\orh(u^{-1})s,\orh(u^{-1})t)-Z_u(\orh(u^{-1})s',\orh(u^{-1})t'))}{2a(|s-s'|^{\alpha/2}+|t-t'|^{\alpha/2})}-1\RT|=0.
\EQN
Further, since for $\alpha/2\in(0,1]$
\BQNY
\abs{\abs{t}^{\alpha/2}-\abs{t'}^{\alpha/2}}\leq \mathbb{C}_1\abs{t-t'}^{\alpha/2},\quad
\abs{\abs{s}^{\alpha/2}-\abs{s'}^{\alpha/2}}\leq \mathbb{C}_2\abs{s-s'}^{\alpha/2}
\EQNY
holds for $t,t'\in[-T_1,T_2],s,s'\in[-S_1,S_2]$,
we have by \eqref{rp1}
\BQNY
&&\sup_{k\in K_u}\underset{(s,t),(s',t')\in D}{\sup_{|(s,t)-(s',t')|<\vn}}M^2_k(u)\E{[Z_u(\orh(u^{-1})s,\orh(u^{-1})t)-Z_u(\orh(u^{-1})s',\orh(u^{-1})t')]Z_u(0,0)}\\
&\quad&\quad \leq\mathbb{C}_3  u^2\underset{(s,t),(s',t')\in D}{\sup_{|(s,t)-(s',t')|<\vn}}\abs{r_{Z_u}(\orh(u^{-1})s,\orh(u^{-1})t)-r_{Z_u}(\orh(u^{-1})s',\orh(u^{-1})t')}\\
&\quad&\quad = \mathbb{C}_3 a  \underset{(s,t),(s',t')\in D}{\sup_{|(s,t)-(s',t')|<\vn}}\abs{\abs{s}^{\alpha/2}+\abs{t}^{\alpha/2}
-\abs{s'}^{\alpha/2}-\abs{t'}^{\alpha/2}}\\
&\quad&\quad \leq \mathbb{C}_3 a  \underset{(s,t),(s',t')\in D}{\sup_{|(s,t)-(s',t')|<\vn}}\LT(\abs{\abs{s}^{\alpha/2}-\abs{s'}^{\alpha/2}}+\abs{\abs{t}^{\alpha/2}
-\abs{t'}^{\alpha/2}}\RT)\\
&\quad&\quad \leq \mathbb{C}_4\vn^{\alpha/2} \rw 0, \ u\rw\IF,\vn\rw 0.
\EQNY
Hence the claim follows from \cite{KEP2015}[Lemma 5.3].
\COM{
Set $D=[-S_1,S_2]\times[-T_1,T_2]$, Conditioning on $\mathcal{A}_{u,k}(w)=\{Z_u(0,0)=M_k(u)-\frac{w}{M_k(u)}\}, w\in \R$, we have for all $u$ large enough
\BQNY
&&\sqrt{2\pi}M_k(u)e^{M_k^2(u)/2}\pk{\sup_{(s,t)\in D(u)}Z_u(s,t)>M_k(u)}\\
&=&\int_{-\IF}^\IF e^{w-\frac{w^2}{2M_k^2(u)}}\pk{\sup_{(s,t)\in D(u)}Z_u(s,t)>M_k(u)\Bigl\lvert \mathcal{A}_{u,k}(w)}dw\\
&=&\int_{-\IF}^\IF e^{w-\frac{w^2}{2M_k^2(u)}}\pk{\sup_{(s,t)\in D}\mathcal{X}_{u,k}(s,t)>w\Bigl\lvert
\mathcal{A}_{u,k}(w)}dw\\
&:=& I_{u,k}
\EQNY
where
\BQNY
\mathcal{X}_{u,k}(s,t)=M_k(u)\LT(Z_u(\orh(u^{-1})s,\orh(u^{-1})t)-M_k(u)\RT)+w.
\EQNY
By the assumption that $\pk{\sup_{(s,t)\in D(u)}Z_u(s,t)>M_k(u)}$ is positive for all $u$ large and any $k\in K_u$, in order to establish the proof we need to show that
\BQN\label{IR2}
\lim_{u\rw\IF}\sup_{k\in K_u}\LT|I_{u,k}-\mathcal{R}_{\alpha}(D)\RT|=0,
\EQN
where $\mathcal{R}_{\alpha}(D)=\E{\sup_{(s,t)\in [-S_1,S_2]\times[-T_1,T_2] }e^{c_0\sqrt{2a}(B_{\frac{\alpha}{2}}(s)+B'_{\frac{\alpha}{2}}(t))-ac_0^2(|s|^{\frac{\alpha}{2}}+|t|^{\frac{\alpha}{2}})}}$.\\
It follows that
\BQNY
&&\sup_{k\in K_u}\LT|I_{u,k}-\mathcal{R}_\alpha(D)\RT|\\
&\leq&\sup_{k\in K_u}\LT|\int_{-W}^W\LT[e^{w-\frac{w^2}{2M^2_k(u)}}\pk{\sup_{(s,t)\in D}\mathcal{X}_{u,k}(s,t)>w\Bigl\lvert\mathcal{A}_{u,k}(w)}
-e^w\pk{\sup_{(s,t)\in D}\zeta(s,t)>w}\RT]dw\RT|\\
&&+\sup_{k\in K_u}\int_{|w|>W}e^{w-\frac{w^2}{2M^2_k(u)}}\pk{\sup_{(s,t)\in D}\mathcal{X}_{u,k}(s,t)>w\Bigl\lvert\mathcal{A}_{u,k}(w)}dw\\
&&+\int_{|w|>W}e^w\pk{\sup_{(s,t)\in D}\zeta(s,t)>w}dw,
\EQNY
for some positive constant $W$, where $\zeta(s,t)=c_0\sqrt{2a}(B_{\frac{\alpha}{2}}(s)+B'_{\frac{\alpha}{2}}(t))-ac_0^2(|s|^{\frac{\alpha}{2}}+|t|^{\frac{\alpha}{2}})$.
Next, we can give the upper bound of each term in the right hand side of the above inequality in the same way as the proof of \nelem{pan2} after we replace (\ref{XM12}), (\ref{XVar}) and (\ref{HC}) as
\BQNY
\E{\mathcal{X}_{u,k}(s,t)\Bigl\lvert\mathcal{A}_{u,k}(w)}\rw -ac^2_0(|s|^{\alpha/2}+|t|^{\alpha/2}),\ u\rw\IF,
\EQNY
\BQNY
\Var\LT(\mathcal{X}_{u,k}(s,t)-\mathcal{X}_{u,k}(s',t')\Bigl\lvert\mathcal{A}_{u,k}(w)\RT)\rw 2ac_0^2(|s-s'|^{\alpha/2}+|t-t'|^{\alpha/2})(1+o(1)),\ u\rw\IF.
\EQNY
\BQNY
\sup_{k\in K_u}\Var\LT(\mathcal{X}_{u,k}(s,t)-\mathcal{X}_{u,k}(s',t')\Bigl\lvert\mathcal{A}_{u,k}(w)\RT)\leq C_1|s-s'|^{\alpha/2}+C_2|t-t'|^{\alpha/2}.
\EQNY
where $C_1$ and $C_2$ are some constants.}
\COM{
Further,
\BQN\label{XM21}
&&\E{\mathcal{X}_{u,k}(s,t)\Bigl\lvert\mathcal{A}_{u,k}(w)}\nonumber\\
&=&-M_k^2(u)(1-r(\orh(u^{-1})s,\orh(u^{-1})t))+w(1-r(\orh(u^{-1})s,\orh(u^{-1})t)),
\EQN
\BQNY
\E{\mathcal{X}_{u,k}(0,0)\Bigl\lvert\mathcal{A}_{u,k}(w)}=\E{\mathcal{X}^2_{u,k}(0,0)\Bigl\lvert\mathcal{A}_{u,k}(w)}=0,
\EQNY
and
\BQNY
&&\Var(\mathcal{X}_{u,k}(s,t)-\mathcal{X}_{u,k}(s',t')\Bigl\lvert\mathcal{A}_{u,k}(w))\\
&=&M_k^2(u)
(\Var(Z_u(\orh(u^{-1})s,\orh(u^{-1})t)-
Z_u(\orh(u^{-1})s',\orh(u^{-1})t'))\\
&&-
(r(\orh(u^{-1})s,\orh(u^{-1})t)-r(\orh(u^{-1})s',\orh(u^{-1})t'))^2).
\EQNY
which does not depend on $w$.\\
From condition (\ref{rp1}) and (\ref{Mku2}), we obtain uniformly with respect to $(s,t)\in D, k\in K_u, w\in[-W,W]$ that
\BQN\label{XM22}
\E{\mathcal{X}_{u,k}(s,t)\Bigl\lvert\mathcal{A}_{u,k}(w)}\rw -ac^2_0(|s|^{\alpha/2}+|t|^{\alpha/2}),\ u\rw\IF,
\EQN
and also for any $(s,t), (s',t')\in D$ uniformly with respect to $k\in K_u$ and any $w\in \mathbb{R}$
\BQN\label{XVar2}
\Var\LT(\mathcal{X}_{u,k}(s,t)-\mathcal{X}_{u,k}(s',t')\Bigl\lvert\mathcal{A}_{u,k}(w)\RT)\rw 2ac_0^2(|s-s'|^{\alpha/2}+|t-t'|^{\alpha/2})(1+o(1)),\ u\rw\IF.
\EQN

Thus, there exists a constant $C_1,C_2>0$ such that for all $(s,t), (s',t')\in D$ and for any $w\in\mathbb{R}$ and all large enough $u$,
\BQN\label{HC2}
\sup_{k\in K_u}\Var\LT(\mathcal{X}_{u,k}(s,t)-\mathcal{X}_{u,k}(s',t')\Bigl\lvert\mathcal{A}_{u,k}(w)\RT)\leq C_1|s-s'|^{\alpha/2}+C_2|t-t'|^{\alpha/2}.
\EQN
Consequently for arbitrary $W>0$,
\BQNY
\lim_{u\rw\IF}\underset{w\in[-W,W]}{\sup_{k\in K_u,}}\LT| \pk{\sup_{(s,t)\in D}\mathcal{X}_{u,k}(s,t)>w\Bigl\lvert\mathcal{A}_{u,k}(w)}
-\pk{\sup_{(s,t)\in D}\zeta(s,t)>w} \RT|=0
\EQNY
and by (\ref{Mku2})
\BQNY
\lim_{u\rw\IF}\sup_{k\in K_u,w\in[-W,W]}e^w[1-e^{-\frac{w^2}{2M^2_k(u)}}]\leq\frac{e^WW^2}{2\lim \inf_{u\rw\IF}\inf_{k\in K_u}M^2_k(u)}\rw 0,\ u\rw\IF,
\EQNY
we obtain
\BQNY
\lim_{u\rw\IF}\sup_{k\in K_u}\LT|\int_{-W}^W\LT[e^{w-\frac{w^2}{2M^2_k(u)}}\pk{\sup_{(s,t)\in D}\mathcal{X}_{u,k}(s,t)>w\Bigl\lvert\mathcal{A}_{u,k}(w)}
-e^w\pk{\sup_{(s,t)\in D}\zeta(s,t)>w}\RT]dw\RT|=0.
\EQNY
Using (\ref{XM21}) and (\ref{XVar2}), for some $\delta\in(0,1)$, $|w|>W$ with $W$ sufficiently large and all $u$ large we have
\BQNY
E_{\mathcal{X}}:=\sup_{k\in K_u,(s,t)\in D}\E{\mathcal{X}_{u,k}(s,t)\Bigl\lvert\mathcal{A}_{u,k}(w)}\leq\delta|w|,
\EQNY
and
\BQNY
\sup_{k\in K_u,(s,t)\in D}\Var\LT(\mathcal{X}_{u,k}(s,t)\Bigl\lvert\mathcal{A}_{u,k}(w)\RT)\leq \mathbb{C}_1
\EQNY
where $\mathbb{C}_1$ is a constant.
Next, by (\ref{HC2}) and Theorem 8.1 of \cite{Pit96}, we obtain for all u large
\BQNY
&&\sup_{k\in K_u}\int_{|w|>W}e^{w-\frac{w^2}{2M^2_k(u)}}\pk{\sup_{t\in[-S_1,S_2]}\mathcal{X}_{u,k}(s,t)>w\Bigl\lvert\mathcal{A}_{u,k}(w)}dw\\
&\leq&\int_{|w|>W}e^{w}\pk{\sup_{t\in[-S_1,S_2]}\mathcal{X}_{u,k}(s,t)-E_{\mathcal{X}}>w-\delta|w|\Bigl\lvert\mathcal{A}_{u,k}(w)}dw\\
&\leq& e^{-W}+\int_W^{\IF}e^w\mathbb{C}_2w^{2/\alpha}\Psi\LT(\frac{w-\delta w}{\mathbb{C}_1}\RT)dw\\
&=:&A_1(W)\rw 0,\ W\rw\IF,
\EQNY
where $\mathbb{C}_2$ is a constant.
Moreover, $$A_2(W):=\int_{|w|>W}e^w\pk{\sup_{t\in[-S_1,S_2]}\zeta(s,t)>w}dw\rw 0,\ W\rw\IF.$$
Hence (\ref{IR2}) follows since
\BQNY
&&\sup_{k\in K_u}\LT|I_{u,k}-\mathcal{R}_\alpha(D)\RT|\\
&\leq&\sup_{k\in K_u}\LT|\int_{-W}^W\LT[e^{w-\frac{w^2}{2M^2_k(u)}}\pk{\sup_{(s,t)\in D}\mathcal{X}_{u,k}(s,t)>w\Bigl\lvert\mathcal{A}_{u,k}(w)}
-e^w\pk{\sup_{(s,t)\in D}\zeta(t)>w}\RT]dw\RT|\\
&&+A_1(W)+A_2(W)\\
&\rw& A_1(W)+A_2(W),\ u\rw\IF,\\
&\rw& 0,\ W\rw\IF.
\EQNY}
\QED

\prooftheo{PreThm1}
We have from {\bf A3}
$$ \lim_{t\rw0}\frac{\rho^2(t)}{t^{2/\lambda}}=
\eta\in[0,\IF],\quad \lim_{u\rw\IF}u^\lambda\orh(u^{-1})=\eta^{-\lambda/2}.$$
%Since $\sigma_u(t_u)=1$ and $\sigma_u(t)$ is continuous function, $\sigma_u(t_u+t)>0$ for $t\in\Delta(u)$ when $u$ large enough.\\
Without loss of \cJ{generality, we consider only the case $t_u= 0$ for $u$ large enough.}\\
By {\bf A2}  for $t\in\Delta(u)$, \kd{for sufficiently large  $u$,}
\BQN\label{sig}
\frac{1}{\mathcal{F}_{u,+\vn}(t)}\leq \sigma_u(t)
\leq\frac{1}{\mathcal{F}_{u,-\vn}(t)},\quad
\mathcal{F}_{u,\pm\vn}(t)=1+u^{-2}\LT[(1\pm\vn)f(u^{\lambda}t)\pm\vn\RT]
\EQN
for small constant $\vn\in(0,1)$. Since further
\BQN\label{pipi}
\pi(u):=\pk{\sup_{t\in \Delta(u)}X_u(t)>M_u}=\pk{\sup_{t\in \Delta(u)}\overline{X_u}(t)\sigma_u(t)>M_u}
\EQN
we have
\BQNY
\pi(u)\leq\pk{\sup_{t\in \Delta(u)}\frac{\overline{X_u}(t)}{\mathcal{F}_{u,-\vn}(t)}>M_u},\ \ \
\pi(u)\geq\pk{\sup_{t\in \Delta(u)}\frac{\overline{X_u}(t)}{\mathcal{F}_{u,+\vn}(t)}>M_u}.
\EQNY
Set for some positive constant $S$
$$I_k(u)=[k\orh(u^{-1})S,(k+1)\orh(u^{-1})S], \ \ k\in\mathbb{Z}.$$
Further, define
\BQNY
&&\mathcal{G}_{u,+\vn}(k)=M_u\sup_{s\in I_k(u)}\mathcal{F}_{u,+\vn}(s),\
N_1(u)=\LT\lfloor \frac{x_1(u)}{S\orh(u^{-1})}\RT\rfloor-\mathbb{I}_{\{x_1\leq0\}},\\
&&\mathcal{G}_{u,-\vn}(k)=M_u\inf_{s\in I_k(u)}\mathcal{F}_{u,-\vn}(s),\
N_2(u)=\LT\lfloor \frac{x_2(u)}{S\orh(u^{-1})}\RT\rfloor+\mathbb{I}_{\{x_2\leq0\}}.
\EQNY
In view of \cite{HP99}, we can find  centered stationary Gaussian processes $Y_{\pm\vn}(t),t\inr$ with continuous trajectories, unit variance and correlation function satisfying
$$
r_{\pm\vn}(t)=1-(1\pm\vn)\rho^2(|t|)(1+o(1)),\ \ t\rw 0.
$$
\underline{Case 1) $\eta=\IF$:}\\
For any $u$ positive
\BQN\label{proof3}
\sum_{k=N_1(u)+1}^{N_2(u)-1}\pk{ \sup_{t\in I_{k}(u)} X_u(t)>M_u}-\sum_{i=1}^2\Lambda_i(u)\leq \pi(u)\leq \sum_{k=N_1(u)}^{N_2(u)}\pk{ \sup_{t\in I_{k}(u)} X_u(t)>M_u},
\EQN
where
$$\Lambda_1(u)=\sum_{k=N_1(u)}^{N_2(u)}\pk{ \sup_{t\in I_{k}(u)} X_u(t)>M_u, \sup_{t\in I_{k+1}(u)} X_u(t)>M_u},$$ and $$\Lambda_2(u)=\sum_{N_1(u)\leq k,l\leq N_2(u), l\geq k+2}\pk{ \sup_{t\in I_{k}(u)} X_u(t)>M_u, \sup_{t\in I_{l}(u)} X_u(t)>M_u}.$$
Set below
$$,\Theta(u)=\frac{\mathcal{H}_{\alpha}}{u^{\lambda}\orh(u^{-1})}\int_{x_1}^{x_2}e^{-f(t)}dt\Psi(M_u).
$$
which is well-defined since  $\int_{x_1}^{x_2}e^{-f(t)}dt<\IF$ follows by the assumption  $f\in C_0^*([x_1,x_2])$.
By Slepian inequality \kd{(see e.g., \cite{AdlerTaylor}),} % for Gaussian random fields, and \cite{GennaSlepian} for stable processes)
\eqref{pipi} and \nelem{pan2}%, we have that
\BQN\label{proof4}
\sum_{k=N_1(u)}^{N_2(u)}\pk{ \sup_{t\in I_{k}(u)} X_u(t)>M_u}
&\leq&\sum_{k=N_1(u)}^{N_2(u)}\pk{ \sup_{t\in I_{k}(u)} \overline{X_u}(t)>\mathcal{G}_{u,-\vn}(k)}\nonumber\\
&\leq&\sum_{k=N_1(u)}^{N_2(u)}\pk{ \sup_{t\in I_{k}(u)} Y_{+\vn}(t)>\mathcal{G}_{u,-\vn}(k)}\nonumber\\
&=&\sum_{k=N_1(u)}^{N_2(u)}\pk{ \sup_{t\in I_{0}(u)} Y_{+\vn}(t)>\mathcal{G}_{u,-\vn}(k)}\nonumber\\
&\sim& \sum_{k=N_1(u)}^{N_2(u)}\mathcal{H}_{\alpha}[0,(1+\vn)^{1/\alpha}S]\Psi(\mathcal{G}_{u,-\vn}(k))\nonumber\\
&\sim& \mathcal{H}_{\alpha}[0,(1+\vn)^{1/\alpha}S]\Psi(M_u)\sum_{k=N_1(u)}^{N_2(u)}e^{-M_{u}^2u^{-2}\inf_{s\in I_k(u)}\LT[(1-\vn)f(u^{\lambda}s)-\vn\RT]}\nonumber\\
&\sim& \frac{\mathcal{H}_{\alpha}[0,(1+\vn)^{1/\alpha}S]}{Su^{\lambda}\orh(u^{-1})}\int_{x_1}^{x_2}e^{- (1-\vn)f(t)+\vn}dt\Psi(M_u)\nonumber\\
&\sim& \Theta(u), \quad  %\frac{\mathcal{H}_{\alpha}}{u^{\lambda}\orh(u^{-1})}\int_{x_1}^{x_2}e^{-f(t)}dt\Psi(M_u), \
u\rw\IF, S\rw\IF,\vn\rw 0.
\EQN
Similarly, we  derive that
\BQN\label{proof5}
\sum_{k=N_1(u)+1}^{N_2(u)-1}\pk{ \sup_{t\in I_{k}(u)} X_u(t)>u}\geq (1+o(1))\Theta(u), %\frac{\mathcal{H}_{\alpha}}{u^{\lambda}\orh(u^{-1})}\int_{x_1}^{x_2}e^{-f(t)}dt\Psi(M_u),\
u\rw\IF,\ S\rw\IF,\ \vn\rw 0.
\EQN
Moreover,
\BQN\label{proof6}
\Lambda_1(u)&\leq&\sum_{k=N_1(u)}^{N_2(u)}\left(\pk{ \sup_{t\in I_{k}(u)} Y_{+\vn}(t)>\widehat{\mathcal{G}}_{u,-\vn}(k)}+\pk{\sup_{t\in I_{k+1}(u)}Y_{+\vn}(t)>\widehat{\mathcal{G}}_{u,-\vn}(k)}\RT.\nonumber\\
& &\LT.-\pk{\sup_{t\in I_{k}(u)\cup I_{k+1}(u)} Y_{-\vn}(t)>\overline{\mathcal{G}}_{u,+\vn}(k)}\right)\nonumber\\
&\leq & \sum_{k=N_1(u)}^{N_2(u)}\left(2\mathcal{H}_{\alpha}[0,(1+\vn)^{1/\alpha}S]-\mathcal{H}_{\alpha}
[0,2(1-\vn)^{1/\alpha}S]\right)\Psi(\widehat{\mathcal{G}}_{u,-\vn}(k))\nonumber\\
&\leq&\left(2\mathcal{H}_{\alpha}[0,(1+\vn)^{1/\alpha}S]-\mathcal{H}_{\alpha}
[0,2(1-\vn)^{1/\alpha}S]\right)\sum_{k=N_1(u)}^{N_2(u)}\Psi(\widehat{\mathcal{G}}_{u,-\vn}(k))\nonumber\\
&=&o(\Theta(u))%o\left(\frac{1}{u^{\lambda}\orh(u^{-1})}\int_{x_1}^{x_2}e^{-f(t)}dt\Psi(M_u)\right)
, \ u\rw\IF, S\rw \IF, \vn\rw 0,
\EQN
where
$$\widehat{\mathcal{G}}_{u,-\vn}(k)=\min (\mathcal{G}_{u,-\vn}(k),\mathcal{G}_{u,-\vn}(k+1)),
\quad \overline{\mathcal{G}}_{u,+\vn}(k)=\max(\mathcal{G}_{u,+\vn}(k),\mathcal{G}_{u,+\vn}(k+1)).$$
By \textbf{A3} for any $(s,t)\in I_{k}(u)\times I_{l}(u)$ with $N_1(u)\leq k,l\leq N_2(u), l\geq k+2$ we have
\BQNY
2\leq Var\left( \overline{X_u}(s)+\overline{X_u}(t)\right)=4-2(1-r_u(s,t))\leq 4-\rho^2(|t-s|)\leq 4-\mathbb{C}_1u^{-2}|(l-k-1)S|^{\alpha/2}
\EQNY
and for $(s,t), (s',t')\in I_{k}(u)\times I_{l}(u)$ with $N_1(u)\leq k,l\leq N_2(u)$
\BQNY
&&1-Cov\LT(\frac{\overline{X_u}(s)+\overline{X_u}(t)}{\sqrt{Var\left( \overline{X_u}(s)+\overline{X_u}(t)\right)}}, \frac{\overline{X_u}(s')+\overline{X_u}(t')}{\sqrt{Var\left( \overline{X_u}(s')+\overline{X_u}(t')\right)}}\RT)\\
&\quad&\quad=\frac{1}{2}\E{\LT(\frac{\overline{X_u}(s)+\overline{X_u}(t)}{\sqrt{Var\left( \overline{X_u}(s)+\overline{X_u}(t)\right)}}-\frac{\overline{X_u}(s')+\overline{X_u}(t')}{\sqrt{Var\left( \overline{X_u}(s')+\overline{X_u}(t')\right)}}\RT)^2}\\
&\quad&\quad=\frac{1}{Var\left( \overline{X_u}(s)+\overline{X_u}(t)\right)}\E{\LT(\overline{X_u}(s)-\overline{X_u}(s')+\overline{X_u}(t)
-\overline{X_u}(t')\RT)^2}\\
&\quad&\quad\quad+Var\left( \overline{X_u}(s')+\overline{X_u}(t')\right)\LT(\frac{1}{\sqrt{Var\left( \overline{X_u}(s)+\overline{X_u}(t)\right)}}-\frac{1}{\sqrt{Var\left( \overline{X_u}(s')+\overline{X_u}(t')\right)}}\RT)^2\\
&\quad&\quad\leq2\E{\LT(\overline{X_u}(s)-\overline{X_u}(s')\RT)^2}+2\E{\LT(\overline{X_u}(t)
-\overline{X_u}(t')\RT)^2}+\E{\LT(\overline{X_u}(s)-\overline{X_u}(s')+\overline{X_u}(t)
-\overline{X_u}(t')\RT)^2}\\
&\quad&\quad \leq 8(1-r_u(s,s')+1-r_u(t,t'))\\
&\quad&\quad = 16u^{-2} \left(\Bigl\lvert\frac{s-s'}{\orh(u^{-1})}\Bigl\lvert^{\alpha/2}+\Bigl\lvert\frac{t-t'}{\orh(u^{-1})}\Bigl\lvert^{\alpha/2}\right).
\EQNY
In view of our assumptions, we can find centered homogeneous Gaussian random fields $Z_u(s,t)$ \cJ{with correlation } %such that its correlation satisfies
$$r_{Z_u}(s,t)=\exp \Biggl(-32 u^{-2} \left(\Biggl\lvert\frac{s}{\orh(u^{-1})}\Bigl\lvert^{\alpha/2}+\Bigl\lvert\frac{t}{\orh(u^{-1})}\Bigl\lvert^{\alpha/2}\right) \Biggr).$$
Slepian inequality, \nelem{pan3} and (\ref{proof4}) imply
\BQN\label{proof7}
\Lambda_2(u)&\leq& \sum_{N_1(u)\leq k,l\leq N_2(u), l\geq k+2}\pk{ \sup_{s\in I_{k}(u)} X_u(s)>M_u, \sup_{t\in I_{l}(u)} X_u(t)>M_u}\nonumber\\
&\leq& \sum_{N_1(u)\leq k,l\leq N_2(u), l\geq k+2}\pk{ \sup_{(s,t)\in I_{k}(u)\times I_{l}(u) } (\overline{X_u}(s)+\overline{X_u}(t))>2\widetilde{\mathcal{G}}_{u,-\vn}(k,l)}\nonumber\\
&\leq& \sum_{N_1(u)\leq k,l\leq N_2(u), l\geq k+2}\pk{ \sup_{(s,t)\in I_{0}(u)\times I_{0}(u) } Z_u(s,t)>\frac{2\widetilde{\mathcal{G}}_{u,-\vn}(k,l)}{\sqrt{4-\mathbb{C}_1u^{-2}|(l-k-1)S|^{\alpha/2}}}}\nonumber\\
&\leq&\sum_{N_1(u)\leq k,l\leq N_2(u), l\geq k+2}\left(\mathcal{H}_{\alpha/2}[0,32^{2/\alpha}S]\right)^2\Psi\left(\frac{2\widetilde{\mathcal{G}}_{u,-\vn}(k,l)}
{\sqrt{4-\mathbb{C}_1u^{-2}|(l-k-1)S|^{\alpha/2}}}\right)\nonumber\\
&\leq& 2\sum_{k=N_1(u)}^{N_2(u)}\sum_{l=1}^{N_2(u)-N_1(u)} \left(\mathcal{H}_{\alpha/2}[0,32^{2/\alpha}S]\right)^2\Psi\left(\frac{2\mathcal{G}_{u,-\vn}(k)}{\sqrt{4-\mathbb{C}_1u^{-2}(l S)^{\alpha/2}}}\right)\nonumber\\
&\leq& 2\sum_{k=N_1(u)}^{N_2(u)} \left(\mathcal{H}_{\alpha/2}[0,32^{2/\alpha}S]\right)^2\Psi\left(\mathcal{G}_{u,-\vn}(k)\right)\sum_{l=1}^\IF e^{-\mathbb{C}_2(lS)^{\alpha/2}}\nonumber\\
&\leq& 2\mathcal{H}_{\alpha/2}32^{2/\alpha} Se^{-\mathbb{C}_3S^{\alpha/2}}\sum_{k=N_1(u)}^{N_2(u)} \mathcal{H}_{\alpha/2}[0,32^{2/\alpha}S]\Psi\left(\mathcal{G}_{u,-\vn}(k)\right)\nonumber\\
&=&o(\Theta(u)), %o\left(\frac{\Psi(M_u)}{u^{\lambda}\orh(u^{-1})}\int_{x_1}^{x_2}e^{-f(t)}dt\right),  \ u\rw\IF,
\quad u\to \IF, S\rw\IF, \vn\rw 0,
\EQN
where $\widetilde{\mathcal{G}}_{u,-\vn}(k,l)=\min (\mathcal{G}_{u,-\vn}(k),\mathcal{G}_{u,-\vn}(l))$.
Combing (\ref{proof3})-(\ref{proof6}) with (\ref{proof7}), we obtain
$$\pi(u)\sim \Theta(u), %\frac{\mathcal{H}_{\alpha}}{u^{\lambda}\orh(u^{-1})}\int_{x_1}^{x_2}e^{-f(t)}dt\Psi(M_u), \ \
\quad u\rw\IF. $$
{\underline{Case 2) $\eta\in(0,\IF)$:} This implies $\lambda=2/\alpha$. \\
	Set for any small constant $\theta\in(0,1)$ and any constant $S_1>0$
	\BQN\label{SS1}
	S_1^*=
	\LT\{
	\begin{array}{ll}
		-S_1,&\ \hbox{if} \ x_1=-\IF;\\
		(x_1+\theta)\eta^{1/\alpha},&\ \hbox{if} \ x_1\in(-\IF,\IF),
	\end{array}
	\RT.\ \
	S_2^*=
	\LT\{
	\begin{array}{ll}
		(x_2-\theta)\eta^{1/\alpha},&\ \hbox{if} \ x_2\in(-\IF,\IF);\\
		S_1,&\ \hbox{if} \ x_2=\IF,
	\end{array}
	\RT.\ \
	\EQN
	\BQN\label{SS2}
	S_1^{**}=
	\LT\{
	\begin{array}{ll}
		-S,&\ \hbox{if} \ x_1=-\IF;\\
		(x_1-\theta)\eta^{1/\alpha},&\ \hbox{if} \ x_1\in(-\IF,\IF),
	\end{array}
	\RT.\ \
	S_2^{**}=
	\LT\{
	\begin{array}{ll}
		(x_2+\theta)\eta^{1/\alpha},&\ \hbox{if} \ x_2\in(-\IF,\IF);\\
		S,&\ \hbox{if} \ x_2=\IF.
	\end{array}
	\RT.\ \
	\EQN
With $K^*=[\orh(u^{-1})S_1^*,\orh(u^{-1})S_2^*]$ and $K^{**}=[\orh(u^{-1})S_1^{**},\orh(u^{-1})S_2^{**}]  $ we have for any $S_1>0$ and $u$ large enough
\BQN
&&\pi(u)\geq\pk{\sup_{t\in K^* } X_u(t)>M_u}\label{proofI0},\\
&&\pi(u)\leq\pk{ \sup_{t\in K^{**}} X_u(t)>M_u}+\underset{k\neq 0,-1}{\sum_{k=N_1(u)}^{N_2(u)}}\pk{ \sup_{t\in I_{k}(u)} X_u(t)>M_u}.\label{proofI1}
\EQN
Using Slepian inequality and Lemma \ref{pan2}, we have that
\BQNY
\pk{ \sup_{t\in K^* } X_u(t)>M_u}&\geq& \pk{ \sup_{t\in K^*} \frac{Y_{-\vn}(t)}{\mathcal{F}_{u,+\vn}(t)}>M_u}\\
%&\sim& \E{\sup_{t\in[S_1^*,S_2^*]}e^{\sqrt{2}c_0B_\alpha(t)-c_0^2((1+\vn)f(\eta^{-1/\alpha} t)+|t|^\alpha)}}\Psi(M_u)\\
&\sim&\mathcal{P}_{\alpha,1}^{h_{+\vn}}[S_1^*,S_2^*]\Psi(M_u),\ u\rw\IF,
\EQNY
where $h_{\pm\vn}(t)=(1\pm\vn)f(\eta^{-1/\alpha}t)\pm \vn$, and similarly
\BQN\label{proofI2}
\pk{ \sup_{t\in K^{**}} X_u(t)>M_u}&\leq& \pk{ \sup_{t\in K^{**}} \frac{Y_{+\vn}(t)}{\mathcal{F}_{u,-\vn}(t)}>M_u}\nonumber\\
&\sim& \mathcal{P}_{\alpha,1}^{h_{-\vn}}[S_1^{**},S_2^{**}]\Psi(M_u),\ u\rw\IF.
\EQN
Moreover, in light of \eqref{eq:fff}, the Slepian inequality and \nelem{pan2}
\BQN\label{proofI3}
\underset{k\neq -1,0}{\sum_{k=N_1(u)}^{N_2(u)}}\pk{ \sup_{t\in I_{k}(u)} X_u(t)>M_u}
&\leq& \underset{k\neq -1,0}{\sum_{k=N_1(u)}^{N_2(u)}} \pk{ \sup_{t\in I_{k}(u)} \frac{Y_{+\vn}(t)}{\mathcal{F}_{u,-\vn}(t)}>M_u}\nonumber\\
&\leq& \underset{k\neq -1,0}{\sum_{k=N_1(u)}^{N_2(u)}} \pk{ \sup_{t\in I_{0}(u)}Y_{+\vn}(t)>\mathcal{G}_{u,-\vn}(k)}\nonumber\\
&\sim& \underset{k\neq -1,0}{\sum_{k=N_1(u)}^{N_2(u)}}\mathcal{H}_\alpha[0,(1+\vn)^{1/\alpha}S]\Psi\left(\mathcal{G}_{u,-\vn}(k)\right)\nonumber\\
&\sim& \mathcal{H}_\alpha[0,(1+\vn)^{1/\alpha}S]\Psi(M_u)\underset{k\neq -1,0}{\sum_{k=N_1(u)}^{N_2(u)}}e^{-\inf_{s\in[k,k+1]}\LT((1-\vn)f(s\eta^{-1/\alpha} S)-\vn\RT)}\nonumber\\
&\sim& \mathbb{C}_4\mathcal{H}_{\alpha}\Psi(M_u) Se^{-\mathbb{C}_5(\eta^{-1/\alpha} S)^{\epsilon_1/2}}e^{\vn}\nonumber\\
&=&o\left(\Psi(M_u)\right),\ u\rw\IF, S\rw\IF, \vn\rw 0.
\EQN
 Letting $\vn\rw0, \ S_1\rw\IF$, $S\rw\IF$, and $\theta\rw 0$ we obtain
\BQNY
\pi(u)\sim\mathcal{P}_{\alpha,\eta}^{f}[x_1,x_2]\Psi(M_u), \quad u\rw\IF.
\EQNY
Next, if we set $x_1(u)=-\LT(\frac{\ln u}{u}\RT)^{\lambda}, x_2(u)=\LT(\frac{\ln u}{u}\RT)^{\lambda}$, then
$$x_1=-\IF,\quad x_2=\IF, \quad S_1^*=-S_1,\quad S_2^*=S_1,\quad S_1^{**}=-S,\quad S_2^{**}=S.$$
Inserting (\ref{proofI2}), (\ref{proofI3}) into (\ref{proofI1}) and letting $\vn\rw0$ leads to
\BQNY
\lim_{u\rw\IF}\frac{\pi(u)}{\Psi(M_u)}\leq \mathcal{P}_{\alpha,\eta}^{f}[-S,S]+ \mathbb{C}_4\mathcal{H}_{\alpha} Se^{-\mathbb{C}_5(\eta^{-1/\alpha} S)^{\epsilon_1/2}}<\IF.
\EQNY
By (\ref{proofI0}), we have
\BQNY
\lim_{u\rw\IF}\frac{\pi(u)}{\Psi(M_u)}\geq\mathcal{P}_{\alpha,\eta}^{f}[-S_1,S_1]>0.
\EQNY
Letting $S_1\rw\IF,S\rw\IF$ we obtain $$\mathcal{P}_{\alpha,\eta}^{f}(-\IF,\IF)\in(0,\IF),\quad \pi(u)\sim\mathcal{P}_{\alpha,\eta}^{f}(-\IF,\IF)\Psi(M_u),\quad u\rw\IF.$$
{\underline{Case 3) $\eta=0$:}
Note that
\begin{align*}
\pi(u)&\leq \pk{\sup_{t\in ((I_{-1}(u)\cup I_0(u))\cap\Delta(u)) }\overline{X_u}(t)\sigma_u(t)>M_u}+
\underset{k\neq -1, 0}{\sum_{k=N_1(u)}^{N_2(u)}}\pk{\sup_{t\in I_k(u) }\overline{X_u}(t)\sigma_u(t)>M_u}=: J_1(u)+J_2(u).
\end{align*}
By \eqref{sig}
\BQN\label{ff}
\frac{1}{\mathcal{F}_{u,+\vn}(t)}\leq\sigma_u(t)\leq \frac{1}{\mathcal{F}_{u,-\vn}(t)}
\leq \frac{1}{1+u^{-2}\inf_{s\in\Delta(u)}[(1-\vn)f(u^{\lambda}s)-\vn]}
\EQN
holds for all $t\in\Delta(u)$. Hence  \nelem{pan2} implies
\BQNY
J_1(u)&\leq& \pk{\sup_{t\in[-\orh(u^{-1})S,\orh(u^{-1})S]}\overline{X_u}(t)>
M_u\LT(1+u^{-2}\inf_{s\in\Delta(u)}[(1-\vn)f(u^{\lambda}s)-\vn]\RT)}\\
&\leq& \pk{\sup_{t\in[-\orh(u^{-1})S,\orh(u^{-1})S]}Y_{+\vn}(t)>M_u\LT(1+u^{-2}\inf_{s\in\Delta(u)}[(1-\vn)f(u^{\lambda}s)-\vn]\RT)}\\
&\sim&\mathcal{H}_\alpha[0,2(1+\vn)^{1/\alpha}S]\Psi \LT(M_u\LT(1+u^{-2}\inf_{s\in\Delta(u)}[(1-\vn)f(u^{\lambda}s)-\vn]\RT)\RT)\\
&\sim&\mathcal{H}_\alpha[0,2(1+\vn)^{1/\alpha}S]\Psi \LT(M_u\RT)e^{-(1-\vn)\omega^*+\vn}\\
&\sim&\Psi \LT(M_u\RT)e^{-\omega^*},\ u\rw\IF, \ S\rw0,\ \vn\rw 0,
\EQNY
where  $\omega^*=\inf_{t\in[x_1,x_2]}f(t)$.
\COM{Further
\BQNY
\Psi \LT(M_u\LT(1+(1-\vn)u^{-2}\inf_{s\in\Delta(u)}[f_1(u^{\lambda_1}s)+f_2(u^{\lambda_2}s)]\RT)\RT)\sim \Psi \LT(M_u\RT){e^{-\omega(x)}}, \ u \rw\IF, \vn\rw 0.
\EQNY}
Furthermore, by \nelem{pan2}, for any $x>0$
\BQN\label{J2}
J_2(u)&\leq&\underset{k\neq -1, 0}{\sum_{k=N_1(u)}^{N_2(u)}} \pk{ \sup_{t\in I_{0}(u)}Y_{+\vn}(t)>\mathcal{G}_{u,-\vn}(k)}\sim\underset{k\neq -1, 0}{\sum_{k=N_1(u)}^{N_2(u)}}\mathcal{H}_\alpha[0,(1+\vn)^{1/\alpha}S]\Psi\left(\mathcal{G}_{u,-\vn}(k)\right)\nonumber\\
&\leq& 2\mathcal{H}_\alpha[0,(1+\vn)^{1/\alpha}S]\Psi(M_u)\sum_{k=1}^{\IF}e^{-(1-2\vn)(kx S)^{\epsilon_1/2}+2\vn}\nonumber\\
&\leq& \mathbb{C}_6\mathcal{H}_{\alpha}\Psi(M_u) Se^{-\mathbb{C}_7(x S)^{\epsilon_1/2}}=o\left(\Psi(M_u)\right),\ u\rw\IF, x\rw\IF, \ S\rw0,
\EQN
hence
 $$\lim_{u\rw\IF}\frac{\pi(u)}{\Psi(M_u)}\leq e^{-\omega^*}, \quad u\rw\IF.$$
Next, since $f\in C_0^*([x_1,x_2])$ there exists $y(u)\in\Delta(u)$ satisfying $$\lim_{u\rw\IF}y(u)u^{\lambda}=y\in\{z\in[x_1,x_2]: f(z)=\omega^*\}.$$
Consequently, in view of \eqref{ff}
\BQNY
\pi(u)
&\geq& \pk{X_u(y(u))>M_u}\\
&\geq&\pk{\overline{X}_u(y(u))>M_u(1+[(1+\vn)f(u^{\lambda} y(u))+\vn]u^{-2})}\\
&=&\Psi\LT(M_u(1+(1+\vn)[f(u^{\lambda} y(u))+\vn]u^{-2})\RT)\\
&\sim&\Psi \LT(M_u\RT){e^{-f(y)}}, \ u \rw\IF,\ \vn\rw 0,
\EQNY
which implies that
\BQNY%\label{ruin2}
\pi(u)\sim \Psi \LT(M_u\RT){e^{-\omega^*}}, \quad u \rw\IF
\EQNY
establishing the proof. \QED

\prooftheo{PreThm2}
%$\Delta(u)=[-\delta_u,\delta_u]$, $\Delta(u)=[-t_u,\delta_u]$ or $\Delta(u)=[-\delta_u, T-t_u]$.
Clearly, for any $u>0$
\BQNY
\pi(u)\le \pk{\sup_{t\in  [0,T]}  X_u(t) >M_u}\le \pi(u)+ \pi_1(u),
\EQNY
where with $D(u):=[0,T]\setminus (t_u+\Delta(u))$,
\BQNY
\pi(u):=\pk{\sup_{t\in  \Delta(u)}  X_u(t_u+t) >M_u},\ \
 \pi_1(u):=\pk{\sup_{t\in D(u)} X_u(t) >M_u}.
\EQNY
 Next, we derive an upper bound  for $\pi_1(u)$ which will finally imply that
\BQN\label{eq:Pi2.2}
\pi_1(u)=o(\pi(u)),\ \ \ u\to\IF.
\EQN
Thus by {\bf A4}, {\bf A5} and Piterbarg inequality (see e.g.,
\cite{Pit96}[Theorem 8.1], \cite{PitH}[Theorem 3] and \cite{KEP2015}[Lemma 5.1])
\BQN \label{eq:Pitqq2}
\pi_1(u)
&=&\pk{\sup_{t\in D(u)} \overline{X}_u(t) \sigma_u(t) >M_u}\nonumber\\
&\leq& \pk{\sup_{t\in D(u)} \overline{X}_u(t)> M_u+\mathbb{C}_1\frac{p(\ln u)^q}{u}}\nonumber\\
&\leq& \mathbb{C}_2T M_u^{2/\varsigma}\Psi\LT(M_u+\mathbb{C}_1\frac{p(\ln u)^q}{u}\RT)\nonumber\\
&=&o\LT(\Psi\LT(M_u\RT)\RT), \quad  u\rw\IF.
\EQN
Since \textbf{A1'} implies \textbf{A1}, by \netheo{PreThm1} and \textbf{A2}, \textbf{A3}, we have
\BQN\label{ppi2}
\pi(u)\sim\Psi\LT(M_u\RT)\left\{
\begin{array}{ll}
\frac{\mathcal{H}_{\alpha}}{u^\lambda\orh(u^{-1})}\int_{x_1}^{x_2} e^{-f(t)}dt,&\hbox{if} \ \ \eta=\IF,\\
 \mathcal{P}_{\alpha,\eta}^{f}[x_1,x_2],& \hbox{if} \ \  \eta\in(0,\IF),\\
1,&  \hbox{if} \ \ \eta=0,
\end{array}
 \right.
\ u\rw\IF,
\EQN
where the result of case $\eta=0$ comes from the fact that $f(t)\geq 0$ for $t\in [x_1,x_2], f(0)=0$ and $0\in [x_1,x_2]$.

Consequently, it follows from \eqref{eq:Pitqq2} and  \eqref{ppi2}  that \eqref{eq:Pi2.2} holds, and thus the proof is complete.
\QED

\proofprop{Thm3}
Without loss of generality we assume that $g_m=g(t_0)=0$.\\
i) We present first the proof for  $t_0\in(0,T)$.
Let $\Delta(u)=[-\delta(u),\delta(u)]$, where $\delta(u)=\LT(\frac{(\ln u)^q }{u}\RT)^{1/\gamma}$ with some large $q>1$.
By \eqref{eq:gtt0} for $u$ large enough and some small $\vn\in(0,1)$
\BQN\label{ccl}
1+\frac{(1-\vn)c\abs{t}^\gamma}{u}\leq\frac{1}{\sigma_u(t+t_0)}:=\frac{u-g(t+t_0)}{u}=1-\frac{g(t+t_0)}{u}\leq 1+\frac{(1+\vn)c\abs{t}^\gamma}{u}
\EQN
holds for all $t\in[-\theta,\theta], \theta>0$.
 It follows that
\BQNY
\Pi(u)\le \pk{\sup_{t\in  [0,T]}  \Xgt >u} \le \Pi(u)+ \Pi_1(u),
\EQNY
 with
\BQNY
 \Pi_1(u):=\pk{\sup_{t\in([0,T] \setminus  [t_0-\theta,t_0+\theta] }  \Xgt >u},
\EQNY
and
\BQNY
 \Pi(u):=\pk{\sup_{t\in [t_0-\theta,t_0+\theta]}\Xgt >u}=\pk{\sup_{t\in [t_0-\theta,t_0+\theta]}X(t)\frac{u}{u-g(t)} >u}.
\EQNY
By \eqref{ccl}, we may further write
\BQN\label{corpp2}
\lim_{u\rightarrow\IF}\sup_{t\in\Delta(u),t\not= 0}\left|\frac{\frac{1}{\sigma_u(t_0+t)}-1}{ cu^{-1}|t|^\gamma }-1\right|=\lim_{u\rightarrow\IF}\sup_{t\in\Delta(u),t\not= 0}\left|\frac{\frac{1}{\sigma_u(t_0+t)}-1}{ c u^{-2}|u^{1/\gamma}t|^\gamma }-1\right|= 0,
\EQN
and
\BQNY
\inf_{t\in[-\theta,\theta]\setminus\Delta(u)}\frac{1}{\sigma_u(t+t_0)}
\geq 1+ \frac{(1-\vn)c (\ln u)^q}{u^2}.
\EQNY
In addition, from \eqref{stationaryR0} we have that
\BQNY
\lim_{u\rightarrow\IF}\underset{t\not=s}{\sup_{  s,t\in\Delta(u)}}\left|\frac{1-r(t_0+t, t_0+s)}{ a|t-s|^\alpha}-1\right|= 0,
\EQNY
and
\BQNY
\sup_{s,t\in[t_0-\theta,t_0+\theta]}\E{X(t)-X(s))^2}\leq \sup_{s,t\in[t_0-\theta,t_0+\theta]} (2-2r(s,t))\leq \mathbb{C}_1 \abs{t-s}^\alpha
\EQNY
hold when $\theta$ \EHc{is} small enough.
Therefore, by \netheo{PreThm2}
\BQNY
\Pi(u)\sim u^{(\frac{2}{\alpha}-\frac{1}{\gamma})_{+}}\Psi\LT( u\RT)
\left\{
\begin{array}{ll}
\mathcal{H}_{\alpha}a^{\frac{1}{\alpha}}\int_{\pa}^\IF e^{ -c\abs{t}^{\gamma}}dt,&\hbox{if} \ \ \alpha<2\gamma,\\
\mathcal{P}_{\alpha,a}^{c|t|^\gamma}[\pa,\IF),&  \hbox{if} \ \ \alpha=2\gamma,\\
1,&  \hbox{if} \ \ \alpha>2\gamma.
\end{array}
 \right.
\EQNY
Moreover, since $g_\theta:=\sup_{t\in[0,T]\setminus[t_0-\theta,t_0+\theta]}g(t)<0$ we have
\BQNY
 \Pi_1(u)
 %&\le &\pk{\sup_{t\in[0,T]}X(t)>u-g(t_0+\delta(u))}\\
 \le \pk{\sup_{t\in[0,T] \setminus [t_0-\theta,t_0+\theta]}X(t)>u-g_\theta  }
 \sim\H_\alpha\ \int_{0}^T \frac{1}{a(t)}dt\ u^{\frac{2}{\alpha}}\Psi\LT(u-g_\theta\RT)=o(\Pi(u)),\ u\rw\IF,
\EQNY
hence the claims follow.\\
For $t_0=0$ and $t_0=T$, we just need to replace $\Delta(u)$ by $\Delta(u)=[0,\delta(u)]$ and $\Delta(u)=[-\delta(u),0]$, respectively. \\
%%%%%
\COM{
	\underline{ii)} \emph{Upper bound}. Set $T_0=0, T_j=\frac{t_j+t_{j+1}}{2},j=1,\ldots,d-1,$ and $T_d=T$, then $t_j\in(T_{j-1},T_{j}),j=1,\ldots,d$.
By i) we have
\BQNY
\pk{\sup_{t\in[T_{j-1},T_j]} \Xgt>u}\sim u^{(\frac{2}{\alpha}-\frac{1}{\gamma_j})_{+}}\Psi\LT(u-g_m\RT)
\LT\{
\begin{array}{ll}
\mathcal{H}_{\alpha}a_j^{\frac{1}{\alpha}}\int_{-\IF}^\IF e^{-c_j|t|^{\gamma_j}}dt,&\ \ \text{if}\ \alpha<2\gamma_j,\\
\mathcal{P}^{c_j|t|^{\gamma_j}}_{\alpha,a_j}(-\IF,\IF), &\ \ \text{if}\ \alpha=2\gamma_j,\\
1,&\ \ \text{if}\ \alpha>2\gamma_j.
\end{array}
\RT.
\EQNY
But by the Bonferroni Inequality we have
\BQNY
\pk{\sup_{t\in  [0,T]}   \Xgt >u}\leq\sum_{j=1}^d\pk{\sup_{t\in[T_{j-1},T_j]} \Xgt>u}.
\EQNY
Thus the upper bound follows.\\
\emph{Lower bound}. Since $r(s,t)<1$ for all $ s\neq t \in[0,T]$, then
$$\theta:=1-\max_{1\le j, k\le d, j\neq k}r(t_j,t_k)>0.$$
Further, %since $r(s,t)$ is continuous function for $(s,t)\in[0,T]\times[0,T]$,
we can find small enough $\vn_j, j=1,2,\ldots,d$ such that
$$
r(s,t)\leq1-\frac{\theta}{2},\ \forall\ (s,t)\in[t_j-\vn_j,t_j+\vn_j]\times[t_k-\vn_k,t_k+\vn_k],$$
and $$[t_j-\vn_j,t_j+\vn_j]\cap[t_k-\vn_k,t_k+\vn_k]=\emptyset,$$
for $ 1\leq j,k\leq d, j\neq k.$ Thus,
\BQNY
\Var(X(s)+X(t))=2+2r(s,t)\leq4-\theta
\EQNY
for $(s,t)\in[t_j-\vn_j,t_j+\vn_j]\times[t_k-\vn_k,t_k+\vn_k],1\leq j,k\leq d, j\neq k.$
\COM{where $\widetilde{a}=\sup_{t\in[0,T]}a(t)$ and $\vn$ is some positive constant.
 Let $Z(s,t)$ be a centered stationary Gaussian field with   unit variance and correlation function
\BQNY
r_{Z}(s,s',t,t')=1-2(1+\vn)\widetilde{a}|s-s'|^{\alpha}-2(1+\vn)\widetilde{a}|t-t'|^{\alpha},
\EQNY
for $(s,t),(s't')\in[0,T]\times[0,T]$.}

By the Borell-TIS inequality,
\BQN \label{eq:BTIS}
&&{\sum_{1\leq j<k\leq d}}\pk{\sup_{t\in[t_j-\vn_j,t_j+\vn_j]} \Xgt>u,\sup_{t\in[t_k-\vn_k,t_k+\vn_k]} \Xgt>u}\nonumber\\
&\quad&\quad\leq{\sum_{1\leq j<k\leq d}}\pk{\sup_{t\in[t_j-\vn_j,t_j+\vn_j]}X(t)>u-g_m,\sup_{t\in[t_k-\vn_k,t_k+\vn_k]}X(t)>u-g_m}\nonumber\\
&\quad&\quad\leq{\sum_{1\leq j<k\leq d}}\pk{\sup_{(s,t)\in[t_j-\vn_j,t_j+\vn_j]\times[t_k-\vn_k,t_k+\vn_k]}(X(s)+X(t))>2u-2g_m}\nonumber\\
&\quad&\quad\le d   \sum_{1\leq j \leq d}  \exp\left(-\frac{(2u-2g_m- \mathbb{Q}(j))^2}{2(4-\theta)}\right)\nonumber\\
&\quad&\quad = o\LT(\Psi\LT(u-g_m\RT)\RT), \ u\rw\IF,
\EQN
\COM{%%%%%%%%%%%%%%%%%%%%%%%%5555
&=&\underset{j\neq k}{\sum_{1\leq j,k\leq d}}\pk{\underset{r(s,t)\neq -1}{\sup_{(s,t)\in[t_j-\vn_j,t_j+\vn_j]\times[t_k-\vn_k,t_k+\vn_k]}}X(s)+X(t)>2u-g(t_j)-g(t_k)}\\
&\leq&\underset{j\neq k}{\sum_{1\leq j,k\leq d}}\pk{\underset{r(s,t)\neq -1}{\sup_{(s,t)\in[t_j-\vn_j,t_j+\vn_j]\times[t_k-\vn_k,t_k+\vn_k]}}
\frac{X(s)+X(t)}{\sqrt{\Var(X(s)+X(t))}}>\frac{2u-g(t_j)-g(t_k)}{\sqrt{4-\theta}}}\\
&\leq&\underset{j\neq k}{\sum_{1\leq j,k\leq d}}\pk{\sup_{(s,t)\in[t_j-\vn_j,t_j+\vn_j]\times[t_k-\vn_k,t_k+\vn_k]}
Z(s,t)>\frac{2u-g(t_j)-g(t_k)}{\sqrt{4-\theta}}}\\
&=&\underset{j\neq k}{\sum_{1\leq j,k\leq d}}\mathbb{Q}(j,k)u^{4/\alpha}\Psi\LT(\frac{2u-g(t_j)-g(t_k)}{\sqrt{4-\theta}}\RT)\\
&\leq&d\sum_{1\leq j\leq d}\mathbb{Q}(j)u^{4/\alpha}\Psi\LT(\frac{2u-2g(t_j)}{\sqrt{4-\theta}}\RT)\\
}%%%%%%%%%%%%%%%%%%%%%%%55
where $ \mathbb{Q}(j)=\E{\sup_{(s,t)\in[t_j-\vn_j,t_j+\vn_j]\times[t_k-\vn_k,t_k+\vn_k]}(X(s)+X(t))}, j=1,\ldots,d, $ are some constants.
But by the Bonferroni Inequality we derive that
\BQNY
\pk{\sup_{t\in  [0,T]}   \Xgt >u}&\geq&\sum_{j=1}^d\pk{\sup_{t\in[t_j-\vn_j,t_j+\vn_j]} \Xgt>u}\\
&&-{\sum_{1\leq j<k\leq d}}\pk{\sup_{t\in[t_j-\vn_j,t_j+\vn_j]} \Xgt>u,\sup_{t\in[t_k-\vn_k,t_k+\vn_k]} \Xgt>u}.
\EQNY
Since further
\BQNY
\pk{\sup_{t\in[t_j-\vn_j,t_j+\vn_j]} \Xgt>u}\geq\pk{X(t_j)+g(t_j)>u}=\Psi\LT(u-g_m\RT), \EQNY
we have in the light of \eqref{eq:BTIS} that
\BQNY
\pk{\sup_{t\in  [0,T]}   \Xgt >u}\geq\sum_{j=1}^d\pk{\sup_{t\in[t_j-\vn_j,t_j+\vn_j]} \Xgt>u}(1+o(1)),\ u\rw\IF.
\EQNY
By combining the upper inequality with
\BQNY
\pk{\sup_{t\in[t_j-\vn_j,t_j+\vn_j]} \Xgt>u}\sim u^{(\frac{2}{\alpha}-\frac{1}{\gamma_j})_{+}}\Psi\LT(u-g_m\RT)
\LT\{
\begin{array}{ll}
a_j^{\frac{1}{\alpha}}\mathcal{H}_{\alpha}\int_{-\IF}^\IF e^{- c_j|t|^{\gamma_j}}dt,&\ \ \text{if}\ \alpha<2\gamma_j,\\
\mathcal{P}^{c_j|t|^{\gamma_j}}_{\alpha,a_j}(-\IF,\IF), &\ \ \text{if}\ \alpha=2\gamma_j,\\
1,&\ \ \text{if}\ \alpha>2\gamma_j,
\end{array}
\RT.
\EQNY
we conclude that the claims follow.\\
}
{ii)} %First we derive the asymptotic of $\pk{\sup_{t\in[A,B]}\Xgt>u}.$
Applying \cite{Pit96}[Theorem 7.1] we obtain
\BQNY
\pk{\sup_{t\in[A,B]} \Xgt >u}=\pk{\sup_{t\in[A,B]} X(t) >u}
\sim\int_{A}^{B}(a(t))^{1/\alpha}dt\mathcal{H}_{\alpha} u^{\frac{2}{\alpha}}\Psi\LT(u\RT).
\EQNY
Set $\Delta_\vn=[A-\vn, B+\vn]\cap [0,T]$ for some $\vn>0$, then we have
\BQNY
\pk{\sup_{t\in[0,T]}(X(t)+g(t))>u}&\geq&\pk{\sup_{t\in[A,B]}(X(t)+g(t))>u},\\
\pk{\sup_{t\in[0,T]}(X(t)+g(t))>u}&\leq&\pk{\sup_{t\in\Delta_\vn}(X(t)+g(t))>u}
+\pk{\sup_{t\in[0,T]\setminus\Delta_\vn}(X(t)+g(t))>u}.
\EQNY
 Since $g$ is a continuous function and $g_\vn:=\sup_{t\in [0,T]\setminus\Delta_\vn}g(t)<0$
\BQNY
\pk{\sup_{t\in [0,T]\setminus\Delta_\vn} (X(t)+g(t))>u}&\leq& \pk{\sup_{t\in[0,T]\setminus\Delta_\vn} X(t) >u-g_\vn}\\
&\leq& \mathbb{C}_2u^{2/\alpha} \Psi(u-g_\vn)=o\LT(u^{2/\alpha} \Psi(u)\RT), \quad u\rw\IF, \vn\rw 0.
\EQNY
Further, we have
\BQNY
\pk{\sup_{t\in\Delta_\vn}(X(t)+g(t))>u}&\leq&\pk{\sup_{t\in\Delta_\vn} X(t)>u}
\cJ{\sim}\int_{A-\vn}^{B+\vn}(a(t))^{\frac{1}{\alpha}}dt\mathcal{H}_\alpha u^{\frac{2}{\alpha}}\Psi(u)\\
&\sim&\int_{A}^{B}(a(t))^{\frac{1}{\alpha}}dt\mathcal{H}_\alpha u^{\frac{2}{\alpha}}\Psi(u) , \quad u\rw\IF, \vn\rw 0.
\EQNY
Hence the claims follow.
\QED

\proofprop{ruintimelocal}
We \kd{give the proof only for $t_0=0$.} In this case, $x\in(0,\IF)$. By definition
\BQNY
\pk{u^{1/\gamma}(\tau_u-t_0)\leq x\big| \tau_u\leq T}=\frac{\pk{\sup_{t\in[0,u^{-1/\gamma}x]} \Xgt>u}}{\pk{\sup_{t\in[0,T]} \Xgt>u}}.
\EQNY
Set $\Delta(u)=[0,u^{-1/\gamma}x]$. For all $u$ large
\BQNY
\pk{\sup_{t\in\Delta(u)} \Xgt>u}=\pk{\sup_{t\in\Delta(u)} X(t)\frac{u}{u-g(t)}>u}.
\EQNY
Denote $X_u(t)=X(t)\frac{u}{u-g(t)} $ and $\sigma_u(t)=\frac{u}{u-g(t)}$. As in the proof of \neprop{Thm3} i), by \netheo{PreThm1} we obtain
\BQNY
\pk{\sup_{t\in\Delta(u)} \Xgt>u}\sim   u^{(\frac{2}{\alpha}-\frac{1}{\gamma})_{+}}\Psi\LT(u\RT)\left\{
\begin{array}{ll}
a^{\frac{1}{\alpha}}\mathcal{H}_{\alpha}\int_{0}^x e^{ -c\abs{t}^{\gamma}}dt,&\hbox{if} \ \ \alpha<2\gamma,\\
\mathcal{P}_{\alpha,a}^{c|t|^\gamma}[0,x],&  \hbox{if} \ \ \alpha=2\gamma,\\
1,&  \hbox{if} \ \ \alpha>2\gamma.
\end{array}
 \right.
\EQNY
Consequently, by \neprop{Thm3} statement i), the results follow.
\QED

\proofprop{MainThm1}
Clearly, for any $u>0$
\BQNY
 \pk{\sup_{t\in[0,T]}  (X(t)+g(t)) >u}
 = \pk{\sup_{t\in[0,T]}  \overline{X}(t)\frac{m_u(t)}{m_u(t_u)} >\frac{u-g(t_u)}{\sigma(t_u)}},
\EQNY
 and\textbf{ A1'} is satisfied.
By the continuity of $\sigma(t)$, $\lim_{u\rw\IF}t_u=t_0$ and $\sigma(t_0)=1$, we have that
for $u$ large enough
$$\sigma(t_u)>0,\ \hbox{and} \ \frac{u-g(t_u)}{\sigma(t_u)}\sim u, \ u\rw\IF.$$
Set next
$$%m_u(t)=\frac{\sigma(t)}{ 1-g(t)/u}, \quad
\cJ{X_u(t)=}\overline{X}(t)\frac{m_u(t)}{m_u(t_u)}, \quad t\in[0,T],  $$
 which has standard deviation function $\sigma_u(t)=\frac{m_u(t_u+t)}{m_u(t_u)} $ and correlation function $r_u(s,t)=r(s,t)$ satisfying assumptions\textbf{ A2}--\textbf{A4}.
 Further, $\overline{X}_u(t)=\overline{X}(t)$ implies\textbf{ A5}.
 Hence  the claims follow from \netheo{PreThm2}. \QED

\proofprop{kA}
For all $u$ large
\BQN\label{rr}
1-r(t_u+t,t_u+s)=\frac{\E{[X(t_u+t)-X(t_u+s)]^2}-[\sigma(t_u+t)-\sigma(t_u+s)]^2}{2\sigma(t_u+t)\sigma(t_u+s)}.
\EQN
Using \kd{that}
\BQNY
\E{[X(t_u+t)-X(t_u+s)]^2}&=&\E{X'^2(t_u+s)}(t-s)^2+o((t-s)^2),\\
\LT[\sigma(t_u+t)-\sigma(t_u+s)\RT]^2&=&\sigma'^2(t_u+t)(t-s)^2+o((t-s)^2),
\EQNY
we have, as $u\to \IF$
\BQNY
1-r(t_u+t,t_u+s)=\frac{\E{X'^2(t_u+t)}-\sigma'^2(t_u+t)}{2\sigma(t_u+t)\sigma(t_u+s)}(t-s)^2+o((t-s)^2).
\EQNY
Since $D(s,t):=\frac{\E{X'^2(t)}-\sigma'^2(t)}{2\sigma(s)\sigma(t)}$ is continuous at   $(t_0,t_0)$, then setting $D=D(t_0,t_0)$ we obtain
$$
\lim_{u\rightarrow\IF}\underset{t\not=s }{\sup_{t\in\Delta(u),s\in\Delta(u)}}\left|\frac{1-r(t_u+t,t_u+s)}{D|t-s|^2}-1\right|=0,
$$
which implies that \textbf{A3} is satisfied.
Next we suppose that $\sigma(t)>\frac{1}{2}$ for any $t\in[0,T]$, since if we set $E_1=\{t\in[0,T]:\sigma(t)\leq \frac{1}{2}\}$, by Borell-TIS inequality
\begin{align*}
\pk{\sup_{t\in E_1}\Xgt >u}\leq \exp\LT(-2\LT(u-\sup_{t\in[0,T]}g(t)-\mathbb{C}_1\RT)^2\RT)
=o\LT(\Psi\LT(\frac{u-g(t_u)}{\sigma(t_u)}\RT)\RT)
\end{align*}
as $u\rw\IF$, where $\mathbb{C}_1=\E{\sup_{t\in[0,T]}X(t)}<0$.
Further by \eqref{rr}
\begin{align*}
\E{( \overline{X}(t)-\overline{X}(s))^2}\leq2-2r(t,s)\leq 4\LT(\sup_{\theta\in[0,T]}\E{X'^2(\theta)}(t-s)^2-\inf_{\theta\in[0,T]}\sigma'^2(\theta)(t-s)^2\RT),
\end{align*}
then \textbf{A5} is satisfied.
 Consequently, the conditions of \neprop{MainThm1} are satisfied and hence the claim follows. \QED

\proofprop{Thm2}
Without loss of generality we assume that $g(t)$ satisfies (\ref{eq:gtt0}) with $g(t_0)=0$.\\
First we present the proof for $t_0\in(0,T)$. Clearly, $m_u$  attains its maximum at the unique point $t_0$. Further, we have
\BQNY
\frac{m_u(t_0)}{m_u(t_0+t)}-1=\frac{1}{\sigma(t_0+t)} (1-\sigma(t_0+t))-\frac{g(t_0+t)}{u\sigma(t_0+t)} .
\EQNY
Consequently, by \eqref{assump-corre} and \eqref{eq:gtt0}
 \BQN \label{eq:Mut0t}
\frac{m_u(t_0)}{m_u(t_0+t)}=1+\LT(b\abs{t}^\beta +\frac{c}{u}\abs{t}^\gamma\RT)(1+o(1)),
\quad t\rw 0
\EQN
 holds for all $u$ large. Further,
 set
$\Delta(u)=[-\delta(u), \delta(u)]$, where
$\delta(u)=\LT(\frac{(\ln u)^q}{u}\RT)^{2/\beta^*}$ for some constant $q>1$ with $\beta^*=\min(\beta,2\gamma)$, and let $f(t)=b|t|^\beta\mathbb{I}_{\{\beta=\beta^*\}}+c|t|^\gamma\mathbb{I}_{\{2\gamma=\beta^*\}}$.
We have
\BQN\label{corpp1}
\lim_{u\rightarrow\IF}\sup_{t\in\Delta(u),t\not= 0}\left|\frac{\LT(\frac{m_u(t_0)}{m_u(t_0+t)}-1\RT)u^2-f(u^{2/{\beta^*}}t)}{f(u^{2/{\beta^*}}t)+\mathbb{I}_{\{\beta\neq 2\gamma\}}}\right|=
 0.
\EQN
By \eqref{assump-corre}
\BQN\label{eq:hc}
\E{(\overline{X}(t)-\overline{X}(s))^2}
=\E{(\overline{X}(t))^2}+\E{(\overline{X}(s))^2}-2\E{\overline{X}(t)\overline{X}(s)}
= 2-2r(s,t)\leq \mathbb{C}_1|t-s|^\alpha
\EQN
holds for $s,t\in[t_0-\theta,t_0+\theta]$, with $\theta>0$ sufficiently small.  By \eqref{eq:Mut0t}, for any $\vn>0$
\BQN\label{eq:mm}
\frac{m_u(t_0)}{m_u(t_0+t)}\geq 1+ \mathbb{C}_2(1-\vn)\frac{(\ln u)^q}{u}
\EQN
holds for all $t\in[-\theta,\theta]\setminus\Delta(u)$. Further
\BQNY
\Pi(u):=\pk{\sup_{t\in [t_0-\theta,t_0+\theta]}  \Xgt >u}\le\pk{\sup_{t\in  [0,T]}  \Xgt >u}\le \Pi(u)+ \Pi_1(u),
\EQNY
 with
\BQNY
 \Pi_1(u):=\pk{\sup_{t\in([0,T] \setminus  [t_0-\theta,t_0+\theta])}  \Xgt >u}.
\EQNY

By\eqref{corpp1}, \eqref{assump-corre}, \eqref{eq:mm}, \eqref{eq:hc} which imply
 \textbf{A2}--\textbf{A5}
and \neprop{MainThm1}, we have
\BQN
\Pi(u)\sim u^{(\frac{2}{\alpha}-\frac{2}{\beta^*})_{+}}\Psi\LT(u\RT)
\LT\{
\begin{array}{ll}
\mathcal{H}_\alpha a^{1/\alpha}\int_{\pa}^\IF e^{-f(t)}dt,&\ \ \text{if}\ \alpha<\beta^*,\\
\mathcal{P}_{\alpha,a}^{f}[\pa,\IF), &\ \ \text{if}\ \alpha=\beta^*,\\
1,&\ \ \text{if}\ \alpha>\beta^*.
\end{array}
\RT.
\EQN
In order to complete the proof it suffices to show that
\BQNY %\label{eq:Pi12}
\Pi_1(u)=o(\Pi(u)).
\EQNY
Since $\sigma_\theta:=\max_{t\in ([0,T] \setminus  [t_0-\theta,t_0+\theta])} \sigma(t)<1$ , by the Borell-TIS inequality we have
\BQNY
 \Pi_1(u)\le \pk{\sup_{t\in([0,T] \setminus  [t_0-\theta,t_0+\theta])} X(t)>u}
 \le \exp\LT(-\frac{(u-\mathbb{C}_3)^2}{2\sigma^2_\theta}\RT)= o(\Pi(u)),
\EQNY
where $\mathbb{C}_3=\E{\sup_{t\in[0,T]} X(t)}<\IF$.\\
\COM{Moreover, in view of \eqref{eq:Mut0t} and \eqref{eq:rst}, by the Piterbarg's inequality we have
\BQNY
\Pi_2(u)
&=&\pk{\sup_{t\in ( [t_0-\theta,t_0+\theta]\setminus (t_0+\Delta(u)))} \Oxt\frac{M_u(t_0)}{M_u(t)} >M_u(t_0)}\nonumber\\
&\leq& C_2 T(M_u(t_0))^{2/\alpha}\Psi\left(M_u(t_0) \LT(1+(1-\vn)\left[\frac{b}{\sigma_0}\left(\frac{(\ln ( u))^q}{u}\right)^{2\beta/{\beta^*}}+cu^{-1}\left(\frac{(\ln ( u))^q}{u}\right)^{2\gamma/{\beta^*}}\right]\RT)\right)\\
 &=& o(\Pi(u))
\EQNY
 for all $u$ large, with some constant $C_2>0$.}
For the cases $t_0=0$ and $t_0=T$, we just need to replace $\Delta(u)$ by $[0,\delta(u)]$ and $[-\delta(u),0]$, respectively. Hence the proof is complete.
\QED

\proofprop{ruintime} i) We shall present the proof only for the case $t_0\in(0,T)$. In this case, $[x_1,x_2]=\R$. By definition, for any $x\in\R$
\BQNY
\pk{u^{\lambda}(\tau_u-t_u)\leq x\big| \tau_u\leq T}=\frac{\pk{\sup_{t\in[0,t_u+u^{-\lambda}x]} \Xgt>u}}{\pk{\sup_{t\in[0,T]} \Xgt>u}}.
\EQNY
For $u>0$ define
$$X_u(t)=\overline{X}(t_u+t)\frac{m_u(t_u+t)}{m_u(t_u)} , \quad \sigma_u(t)=\frac{m_u(t_u+t)}{m_u(t_u)}. $$
 As in the proof of \neprop{MainThm1}, we obtain
\BQNY
\pk{\sup_{t\in[0,t_u+u^{-\lambda}x]} \Xgt>u}=\pk{\sup_{t\in[0,t_u+u^{-\lambda}x]}X_u(t)>\frac{u-g(t_u)}{\sigma(t_u)}},
\EQNY
and \textbf{A1'}, \textbf{A2}--\textbf{A5} are satisfied with $\Delta(u)=[-\delta_u,u^{-\lambda}x]$.
Clearly, for any $u>0$
\BQNY
\pi(u)\leq\pk{\sup_{t\in[0,t_u+u^{-\lambda}x]}X_u(t)>\frac{u-g(t_u)}{\sigma(t_u)}}
\leq \pi(u)+\pi_1(u),
\EQNY
where
\BQNY
\pi(u)=\pk{\sup_{t\in[t_u-\delta(u),t_u+u^{-\lambda}x]}X_u(t)>\frac{u-g(t_u)}{\sigma(t_u)}},\quad
\pi_1(u)=\pk{\sup_{t\in[0,t_u-\delta(u)]}X_u(t)>\frac{u-g(t_u)}{\sigma(t_u)}}.
\EQNY
Applying  \netheo{PreThm1} we have
\BQN\label{ruin1}
\pi(u)\sim  \Psi\LT(\frac{u-g(t_u)}{\sigma(t_u)}\RT)\left\{
\begin{array}{ll}
	\frac{\mathcal{H}_{\alpha}}{u^\lambda\orh(u^{-1})}\int_{-\IF}^{x} e^{-f(t)}dt,&\hbox{if} \ \ \eta=\IF,\\
	\mathcal{P}_{\alpha,\eta}^{f}(-\IF,x],& \hbox{if} \ \  \eta\in(0,\IF),\\
	\sup_{t\in(-\IF,x]}e^{-f(t)},& \hbox{if} \ \  \eta=0.
\end{array}
\right.
\EQN
 In view of \eqref{eq:Pitqq2}
\BQNY
\pi_1(u)=o\LT(\Psi\LT(\frac{u-g(t_u)}{\sigma(t_u)}\RT)\RT), \quad  u\rw\IF,
\EQNY
hence
\BQNY
\pk{\sup_{t\in[0,t_u+u^{-\lambda}x]} \Xgt>u} \sim  \pi(u), \quad u\rw\IF
\EQNY
and thus the claim follows by \eqref{ruin1} and \neprop{MainThm1}. \\
ii) We give the proof of $t_0=T$. In this case $x\in(-\IF,0)$ implying
\BQNY
\pk{u^{2/\beta^*}(\tau_u-T)\leq x\big| \tau_u\leq T}=\frac{\pk{\sup_{t\in[0,T+u^{-2/\beta^*}x]} \Xgt>u}}{\pk{\sup_{t\in[0,T]} \Xgt>u}}.
\EQNY
Set $\delta_u=\LT(\frac{(\ln u)^q}{u}\RT)^{2/\beta^*}$ for some $q>1$ and let  $$\Delta(u)=[-\delta_u,u^{-2/\beta^*}x], \quad \sigma_u(t)=\frac{m_u(t)}{m_u(T)}, $$ with $$m_u(t)=\frac{\sigma(t)}{1-g(t)/u}, \quad \quad X_u(t)=\overline{X}(t)\frac{m_u(t)}{m_u(T)} .$$
For all $u$ large, we have
\BQNY
\pi(u)\leq\pk{\sup_{t\in[0,T+u^{-2/\beta^*}x]} \Xgt>u}\leq\pi(u)+\pk{\sup_{t\in[0,T-\delta_u]} \Xgt>u},
\EQNY
where
$$\pi(u):=\pk{\sup_{t\in\Delta(u)} (X(T+t)+g(T+t))>u}
=\pk{\sup_{t\in\Delta(u)} X_u(T+t)>u}.$$
As in the proof of \neprop{Thm2} it follows that the Assumptions \textbf{A2}--\textbf{A5} hold with $\Delta(u)=[-\delta_u,u^{-2/\beta^*}x]$. Hence an application of  \netheo{PreThm1} yields
\begin{align}\label{ruintimep1}
\pi(u)
\sim  u^{(\frac{2}{\alpha}-\frac{2}{\beta^*})_{+}}\Psi\LT(u\RT)
\LT\{
\begin{array}{ll}
a^{1/\alpha}\mathcal{H}_\alpha\int_{-x}^{\IF} e^{-f(t)}dt,& \text{if}\ \alpha<\beta^*,\\
\mathcal{P}_{\alpha,a}^{f}[-x,\IF), &  \text{if}\ \alpha=\beta^*,\\
e^{-f(x)},&  \text{if}\ \alpha>\beta^*.
\end{array}
\RT.
\end{align}
In view of \eqref{eq:Pitqq2}
\BQNY
\pk{\sup_{t\in[0,T-\delta_u]} \Xgt>u}=\pk{\sup_{t\in[0,T-\delta_u]} X_u(t)>u}=o\LT(\Psi\LT(u\RT)\RT), \quad  u\rw\IF
\EQNY
implying
\BQNY
\pk{\sup_{t\in[0,T+u^{-2/\beta^*}x]} \Xgt>u}\sim \pi(u), \quad u\rw\IF.
\EQNY
Consequently, the proof follows by \eqref{ruintimep1} and \neprop{Thm2}.
\QED

\proofprop{exmodel}
Set next  $A(t)= \int_0^t e^{-\delta v}dB(v)$ and define
$$\widetilde{U}(t)=u+c \int_0^t e^{-\delta v}dv-\sigma A(t),\ \ t\ge0.$$
Since
$$\sup_{t\in[0,\IF)}\E{[A(t)]^2}=\frac{1}{2\delta}$$
implying $ \sup_{t\in[0,\IF)}\E{| A(t)|}<\IF,$ then  by the martingale convergence theorem in \cite{PP1966} we have that  $\widetilde{U}(\IF):=\lim_{t\rw\IF}\widetilde{U}(t)$ exists and is finite almost surely. Clearly, for any $u>0$
\BQNY
p(u)&=&\pk{\inf_{t\in[0,\IF)}\widetilde{U}(t)<0}\\
%&=&\pk{\inf_{t\in[0,\IF]}\widetilde{U}(t)<0}\\
&=&\pk{\sup_{ t\in[0,\IF]}\LT(\sigma A(t)-c\int_0^t e^{-\delta v}dv\RT)>u}\\
&=& \mathbb{P}\left\{\sup_{t\in[0,1]}\LT(\sigma A(-\frac{1}{2\delta}\ln t) %\int_0^{-\frac{1}{2\delta}\ln t^*} e^{-\delta v}dB(v)
-\frac{c}{\delta}(1-{t}^{\frac{1}{2}})\RT)>u\right\}.
\EQNY
 The proof will follow by applying \neprop{MainThm1}, hence we check next the assumptions therein for this specific model.

Below, we set $Z(t)=\sigma A(-\frac{1}{2\delta}\ln t) $
with variance function given by
\BQNY
V_Z^2(t)=Var\left(\sigma\int_0^{-\frac{1}{2\delta}\ln t} e^{-\delta v}dB(v)\right)=\frac{\sigma^2}{2\delta}(1-t), \quad t\in [0,1].
\EQNY
We show next that for $u$ sufficiently large, the function
$$M_u(t):=\frac{uV_Z(t)}{G_u(t)}=\frac{\frac{\sigma}{\sqrt{2\delta}}\sqrt{1-t}}{1+\frac{c}{\delta u}(1-t^{1/2})},  \quad 0\leq t \leq 1,$$
with $G_u(t):=u+\frac{c}{\delta}(1-t^{\frac{1}{2}})$ attains its maximum at the unique point
$t_u=\LT(\frac{c}{\delta u+c}\RT)^2.$
%We calculate the derivative of $M_u(t)$, i.e,
In fact, we have
\BQN\label{VIF1}
[M_u(t)]_t:=\frac{d M_u(t)}{dt}&=&\frac{dV_Z(t)}{dt}\cdot\frac{u}{G_u(t)}-\frac{V_Z(t)}{G_u^2(t)}
\LT(-\frac{cu}{2\delta}t^{-\frac{1}{2}}\RT)=\frac{u}{2G_u^2(t) V_z(t)}\left[\frac{dV_Z^2(t)}{dt}G_u(t)+V_Z^2(t)\frac{ct^{-\frac{1}{2}}}{\delta}\right]\nonumber\\
&=&\frac{u\sigma^2 t^{-1/2}}{4\delta G_u^2(t) V_Z(t)}\LT[\frac{c}{\delta}-\LT(u+\frac{c}{\delta}\RT)t^{\frac{1}{2}}\RT].
\EQN
Letting $[M_u(t)]_t=0$, we get $t_u=\LT(\frac{c}{\delta u+c}\RT)^2$.
By (\ref{VIF1}), $[M_u(t)]_t>0$ for $t\in(0,t_u)$ and $[M_u(t)]_t<0$ for $t\in(t_u,1]$, so $t_u$ is the unique maximum point of $M_u(t)$ over $[0,1]$. Further
\BQNY
M_u:=M_u(t_u)=\frac{\sigma u}{\sqrt{2\delta u^2+4cu}}=\frac{\sigma}{\sqrt{2\delta}}(1+o(1)),\ u\rw\IF.
\EQNY
We set $\delta(u)=\LT(\frac{(\ln u)^q}{u}\RT)^2$ for some $q>1$, and
$\Delta(u)=[-t_u,\delta(u)].$
Next we check the assumption {\bf A2}. It follows that
\BQNY
\frac{M_u}{M_u(t_u+t)}-1=\frac{[G_u(t_u+t)V_Z(t_u)]^2-[G_u(t_u)V_Z(t_u+t)]^2}{V_Z(t_u+t)G_u(t_u)[G_u(t_u+t)V_Z(t_u)+V_Z(t_u+t)G_u(t_u)]}.
\EQNY
We further write
\BQNY
&&[G_u(t_u+t)V_Z(t_u)]^2-[G_u(t_u)V_Z(t_u+t)]^2\\
&&\quad\quad=\LT[\LT(u+\frac{c}{\delta}\RT)-\frac{c}{\delta}\sqrt{t_u+t}\RT]^2\frac{\sigma^2}{2\delta}(1-t_u)
-\LT[\LT(u+\frac{c}{\delta}\RT)-\frac{c}{\delta}\sqrt{t_u}\RT]^2\frac{\sigma^2}{2\delta}(1-t_u-t)\\
&&\quad\quad=\LT(u+\frac{c}{\delta}\RT)^2\frac{\sigma^2}{2\delta}t-2\LT(u+\frac{c}{\delta}\RT)\frac{c\sigma^2}{2\delta^2}(\sqrt{t_u+t}-\sqrt{t_u})(1-t_u)-\frac{c^2\sigma^2}{2\delta^3}t\\
&&\quad\quad=\LT(u+\frac{c}{\delta}\RT)^2\frac{\sigma^2}{2\delta}t(1-t_u)-2\LT(u+\frac{c}{\delta}\RT)^2\frac{\sigma^2}{2\delta}(1-t_u)\sqrt{t_u}(\sqrt{t_u+t}-\sqrt{t_u})\\
&&\quad\quad=\frac{\sigma^2}{2\delta}\LT[\LT(u+\frac{c}{\delta}\RT)^2-\LT(\frac{c}{\delta}\RT)^2\RT](\sqrt{t+t_u}-\sqrt{t_u})^2\\
&&\quad\quad=\frac{\sigma^2}{2\delta}\LT( u^2+\frac{2c}{\delta}u\RT)(\sqrt{t+t_u}-\sqrt{t_u})^2.
\EQNY
Since for any $t\in\Delta(u)$
\BQNY
\sqrt{\frac{\sigma^2}{2\delta}(1-t_u-\delta(u))}\leq V_Z(t_u+t)\leq\sqrt{\frac{\sigma^2}{2\delta}},\
u+\frac{c}{\delta}-\frac{c}{\delta}\sqrt{t_u+\delta(u)}\leq G_u(t_u+t)\leq u+\frac{c}{\delta},
\EQNY
we have for all large $u$
\BQNY
V_Z(t_u+t)G_u(t_u)[G_u(t_u+t)V_Z(t_u)+V_Z(t_u+t)G_u(t_u)] \leq \frac{\sigma^2}{\delta}\LT(u+\frac{c}{\delta}\RT)^2
\EQNY
and
\BQNY
V_Z(t_u+t)G_u(t_u)[G_u(t_u+t)V_Z(t_u)+V_Z(t_u+t)G_u(t_u)]
&\geq&\frac{\sigma^2}{\delta}(1-t_u-\delta(u))\LT(u+\frac{c}{\delta}-\frac{c}{\delta}\sqrt{t_u+\delta(u)}\RT)^2\\
&\geq& \frac{\sigma^2}{\delta}\LT[\LT(u+\frac{c}{\delta}\RT)^2-u\RT].
\EQNY
Thus as $u\rw\IF$
\BQN\label{MM1}
\inf_{t\in\Delta(u),t\neq 0}\frac{M_u/M_u(t_u+t)-1}{\frac{1}{2}\LT(\sqrt{u^2t+\frac{c^2}{\delta^2}}-\frac{c}{\delta}\RT)^2 u^{-2}}-1
\geq\frac{\frac{1}{2} \frac{u^2+\frac{2c}{\delta}u}{(u+\frac{c}{\delta})^2}(\sqrt{t+t_u}-\sqrt{t_u})^2}{\frac{1}{2}\LT(\sqrt{t+\frac{c^2}{(\delta u)^2}}-\frac{c}{\delta u}\RT)^2 }-1
\geq\frac{u^2+\frac{2c}{\delta}u}{(u+\frac{c}{\delta})^2}-1 \rw 0,
\EQN
where we used the fact that for $t\in\Delta(u)$
$$(\sqrt{t+t_u}-\sqrt{t_u})^2\geq\LT(\sqrt{t+\frac{c^2}{(\delta u)^2}}-\frac{c}{\delta u}\RT)^2.$$ Furthermore, since
\BQNY
0&\leq&\frac{\sqrt{t+t_u}-\sqrt{t_u}}{\sqrt{t+\frac{c^2}{(\delta u)^2}}-\frac{c}{\delta u}}-1
=\frac{\sqrt{t+\frac{c^2}{(\delta u)^2}}+\frac{c}{\delta u}}{\sqrt{t+t_u}+\sqrt{t_u}}-1
\leq\frac{\sqrt{t+\frac{c^2}{(\delta u)^2}}-\sqrt{t+t_u}}{\sqrt{t+t_u}+\sqrt{t_u}}\\
&=&\frac{\frac{c^2}{(\delta u)^2}-t_u}{(\sqrt{t+t_u}+\sqrt{t_u})(\sqrt{t+\frac{c^2}{(\delta u)^2}}+\sqrt{t+t_u})}
\leq\frac{\sqrt{\frac{c^2}{(\delta u)^2}-t_u}}{\sqrt{t_u}}
=\sqrt{\LT(1+\frac{c}{\delta u}\RT)^2-1},
\EQNY
we have as $u\rw\IF$
\BQN\label{MM2}
\sup_{t\in\Delta(u),t\neq 0}\frac{M_u/M_u(t_u+t)-1}{\frac{1}{2}\LT(\sqrt{u^2t+\frac{c^2}{\delta^2}}-\frac{c}{\delta}\RT)^2 u^{-2}}-1
&\leq&\frac{\frac{1}{2} \frac{u^2+\frac{2c}{\delta}u}{(u+\frac{c}{\delta})^2-u}(\sqrt{t+t_u}-\sqrt{t_u})^2}{\frac{1}{2}\LT(\sqrt{t+\frac{c^2}{(\delta u)^2}}-\frac{c}{\delta u}\RT)^2 }-1\nonumber\\
&\leq&\frac{u^2+\frac{2c}{\delta}u}{(u+\frac{c}{\delta})^2-u}\LT(1+\sqrt{\LT(1+\frac{c}{\delta u}\RT)^2-1}\RT)^2-1\rw 0.
\EQN
Consequently, (\ref{MM1}) and (\ref{MM2}) imply
\BQN\label{eq:MM12}
\lim_{u\rw\IF}\sup_{t\in\Delta(u),t\neq 0}\LT|\frac{M_u/M_u(t_u+t)-1}{\frac{1}{2}\LT(\sqrt{u^2t+\frac{c^2}{\delta^2}}-\frac{c}{\delta}\RT)^2 u^{-2}}-1\RT|=0.%\ \ \    \textbf{(Assumption A2)}
\EQN
Since for $0\leq t'\leq t<1$, the correlation function of $Z(t)$ equals
\BQNY
r(t,t')=\frac{\E{(\sigma\int_0^{-\frac{1}{2\delta}\ln t} e^{-\delta v}dB(v))(\sigma\int_0^{-\frac{1}{2\delta}\ln t'} e^{-\delta v}dB(v))}}{\sqrt{\frac{\sigma^2}{2\delta}(1-t)}\sqrt{\frac{\sigma^2}{2\delta}(1-t')}}
=\frac{\sqrt{1-t}}{\sqrt{1-t'}}=1-\frac{t-t'}{\sqrt{1-t'}(\sqrt{1-t'}+\sqrt{1-t})},
\EQNY
we have
\BQN\label{ar}
\sup_{t,t'\in\Delta(u),t'\neq t}\LT|\frac{1-r(t_u+t,t_u+t')}{\frac{1}{2}|t-t'|}-1\RT|
&=&\sup_{t,t'\in\Delta(u),t'\neq t}\LT|\frac{2}{\sqrt{1-t-t_u}(\sqrt{1-t'-t_u}+\sqrt{1-t-t_u})}-1\RT|\nonumber\\
&\leq&\frac{1}{1-(\frac{c}{c+\delta u})^2-(\frac{(\ln u)^q}{u})^2}-1 \rw 0,\ \ u\rw\IF.%\ \ \textbf{(Assumption A3)}
\EQN
Further, for some small $\theta\in(0,1)$, we obtain (set below $\overline{Z}(t)=\frac{Z(t)}{V_Z(t)}$)
\BQN\label{hh}
\mathbb{E}\left(\overline{Z}(t)-\overline{Z}(t')\right)^2=2-2r(t,t')\leq\mathbb{C}_1|t-t'|%\ \  \ \textbf{(Assumption A5)}
\EQN
for $ t, t'\in[0,\theta]$.
\def\TPi{\widetilde{\Pi}}
 For all $u$ large %Set $\Delta_\theta=[0,\theta_1]$, then
\BQNY
\Pi(u):=\pk{\sup_{t\in[0,\theta]}\LT(Z(t)-\frac{c}{\delta}(1-t^{\frac{1}{2}})\RT)>u}\leq p(u)\leq\Pi(u) +\TPi(u),
\EQNY
where
\BQNY
\TPi(u):=\pk{\sup_{t\in [\theta,1]}\LT(Z(t)-\frac{c}{\delta}(1-t^{\frac{1}{2}})\RT)>u}\leq\pk{\sup_{t\in [\theta,1]}Z(t)>u}.
\EQNY
%where $Z_u(t):=\frac{uZ(t)}{G_u(t)M_u}$.
%In the following, we show that
%shall focus on the asymptotic of $\Pi(u)$ as $u\to\IF$, and finally show that
Moreover, for all $u$ large
\BQN\label{assA41}
\frac{1}{M_u(t)}-\frac{1}{M_u}&\geq& \frac{[G_u(t)V_Z(t_u)]^2-[G_u(t_u)V_Z(t)]^2}{2uV_Z^3(t_u)G_u(t_u)}
=\frac{\frac{\sigma^2}{2\delta}(u^2+\frac{2c}{\delta}u)(\sqrt{t}-\sqrt{t_u})^2}{2u[\frac{\sigma^2}{2\delta}(1-t_u)]^{3/2}[u+\frac{c}{\delta}(1-\sqrt{t_u})]}\nonumber\\
%&=&\frac{\sqrt{2}(\delta u+c)^4}{2u\delta^2\sigma(\delta u^2+2cu)^{3/2}}\frac{t^2}{(\sqrt{t+t_u}+\sqrt{t_u})^2}\nonumber\\
&\geq&\mathbb{C}_2 (\sqrt{t}-\sqrt{t_u})^2
\geq\frac{\mathbb{C}_2\delta^2(u)}{\LT(\sqrt{\delta(u)+t_u}+\sqrt{t_u}\RT)^2}
\geq\mathbb{C}_3\frac{(\ln u)^{2q}}{u^2}%\ \ \textbf{(Assumption A4)}
\EQN
holds for any $t\in\LT[t_u+\delta(u),\theta\RT]$, therefore
\BQNY
\inf_{t\in\LT[t_u+\delta(u),\theta\RT]}\frac{M_u}{M_u(t)}\geq 1+\mathbb{C}_3\frac{(\ln u)^{q}}{u^2}.
\EQNY
The above inequality combined with (\ref{eq:MM12}), (\ref{ar}), (\ref{hh}) and \neprop{MainThm1} yields
\BQNY
\Pi(u)\sim\mathcal{P}_{1,\delta/\sigma^2}^{h}\LT[-\frac{c^2}{\delta^2},\IF\RT)\Psi\LT(\frac{1}{\sigma}\sqrt{2\delta u^2+4cu}\RT), \ u\rw\IF.
\EQNY
Finally, since
\BQNY
\sup_{t\in [\theta,1]}V^2_Z(t)
%=\sup_{t\in [\theta,1]}\frac{\sigma^2u^2(1-t)}{2\delta M_u^2( u+\frac{c}{\delta}(1-t^{1/2}))^2}
\leq \frac{\sigma^2}{2\delta}(1-\theta),\ \ \hbox{and} \ \
\E{\sup_{t\in [\theta,1]}Z(t)}\leq \mathbb{C}_4<\IF,
\EQNY
by Borell-TIS inequality
\BQNY
\TPi(u)\leq\pk{\sup_{t\in [\theta,1]}Z(t)>u}
\leq \exp\LT(-\frac{\delta(u-\mathbb{C}_4)^2}{\sigma^2(1-\theta)}\RT)
%=o\LT(\Psi\LT(\frac{1}{\sigma}\sqrt{2\delta u^2+4cu}\RT)\RT)
=o(\Pi(u)),\quad u\rw \IF,
\EQNY
which establishes the proof.
Next, we consider that
\begin{align*}
\pk{u^{2}\LT(e^{-2\delta\tau_u}-\LT(\frac{c}{\delta u+c}\RT)^2\RT)\leq x\Big| \tau_u<\IF}
&=\frac{\pk{\inf_{t\in[-\frac{1}{2\delta}\ln \LT(t_u+u^{-2}x\RT),\IF)}\widetilde{U}(t)<0}}
{\pk{\inf_{t\in[0,\IF)}\widetilde{U}(t)<0}}\\
&=\frac{\mathbb{P}\left\{\sup_{t\in[0,t_u+u^{-2}x]}\LT(\sigma A(-\frac{1}{2\delta}\ln t)%\int_0^{-\frac{1}{2\delta}\ln t^*} e^{-\delta v}dB(v)
	-\frac{c}{\delta}(1-{t}^{\frac{1}{2}})\RT)>u\right\}}
{\mathbb{P}\left\{\sup_{t\in[0,1]}\LT(\sigma A(-\frac{1}{2\delta}\ln t) %\int_0^{-\frac{1}{2\delta}\ln t^*} e^{-\delta v}dB(v)
	-\frac{c}{\delta}(1-{t}^{\frac{1}{2}})\RT)>u\right\}}\\
&=\pk{u^{2}\LT(\tau^*_u-t_u\RT)\leq x\big| \tau^*_u<1},
\end{align*}
where
$$\tau^*_u=\{t \in [0,1]: \sigma A( -\frac{1}{2\delta}\ln t)-\frac{c}{\delta}(1-{t}^{\frac{1}{2}})>u\}.$$
 The proof follows by \neprop{ruintime} i). \QED

\section{Appendix}
{\bf Proof of \eqref{prop1}:}
Let $\xi(t), t\in\R$ be a centered stationary Gaussian process with   unit variance and correlation function $r$ satisfying
\begin{align*}
1-r(t)\sim a|t|^\alpha, \ t\rw0,\ \ a>0, \ \alpha\in(0,2].
\end{align*}
In view of by \netheo{PreThm1}, for $-\IF<x_1< x_2<\IF$ and $f\in C^*_0([x_1,x_2])$ we have
\begin{align*}
\pk{\sup_{t\in[u^{-2/\alpha}x_1,u^{-2/\alpha}x_2]}\frac{\xi(t)}{1+u^{-2}f(u^{2/\alpha}t)}>u}\sim\Psi(u)\mathcal{P}_{\alpha,a}^{f}[x_1,x_2], \quad u\to\IF
\end{align*}
and for any $y\in\R$
\begin{align*}
&\pk{\sup_{t\in[u^{-2/\alpha}x_1,u^{-2/\alpha}x_2]}\frac{\xi(t)}{1+u^{-2}f(u^{2/\alpha}t)}>u}\\
&=\pk{\sup_{t\in[u^{-2/\alpha}(x_1-y),u^{-2/\alpha}(x_2-y)]}
	\frac{\xi(t+yu^{-2/\alpha})(1+u^{-2}f(y))}{1+u^{-2}f(y+u^{2/\alpha}t)}>u(1+u^{-2}f(y))}\\
&\sim \Psi(u(1+u^{-2}f(y)))
\mathcal{P}_{\alpha,a}^{f_y(t)-f(y)}[x_1-y,x_2-y]\\
&\sim \Psi(u)
\mathcal{P}_{\alpha,a}^{f_y(t)}[x_1-y,x_2-y].
\end{align*}
\kd{Let}
$$Z_u(t)=\frac{\xi(t+yu^{-2/\alpha})(1+u^{-2}f(y))}{1+u^{-2}f(y+u^{2/\alpha}t)},\quad   t\in[u^{-2/\alpha}(x_1-y),u^{-2/\alpha}(x_2-y)]$$
and denote its variance function by $\sigma^2_{Z_u}(t)$. \kd{Then}
\begin{align*}
\LT(\frac{1}{\sigma_{Z_u}(t)}-1\RT)u^2=
\LT(\frac{1+u^{-2}f(y+u^{2/\alpha}t)}{1+u^{-2}f(y)}-1\RT)u^2=\frac{f(y+u^{2/\alpha}t)-f(y)}
{1+u^{-2}f(y)},
\end{align*}
i.e.,
\BQNY
\lim_{u\rw\IF}\sup_{t\in[u^{-2/\alpha}(x_1-y),u^{-2/\alpha}(x_2-y)]}
\abs{\frac{\LT(\frac{1}{\sigma_{Z_u}(t)}-1\RT)u^2}{f(y+u^{2/\alpha}t)-f(y)}-1}=0.
\EQNY
Consequently, we have
$$\mathcal{P}_{\alpha,a}^{f}[x_1,x_2]=\mathcal{P}_{\alpha,a}^{f_y}[x_1-y,x_2-y].$$
Further,  letting $x_2\rw \IF$ yields  $\mathcal{P}_{\alpha,a}^{f}[x_1,\IF)=\mathcal{P}_{\alpha,a}^{f_y}[x_1-y,\IF)$.
\kd{This completes the proof.} \QED

{\bf Proof of Example \ref{exam2}:}
We have $t_0=0, \gamma=1,g_m=0$.
Then by \neprop{Thm3} statement i)
\BQNY
\pk{\max_{t\in[0,T]}\LT(X(t)-ct\RT)>u }\sim\Psi(u )
\LT\{\begin{array}{ll}
c^{-1}a^{1/\alpha} u^{2/\alpha-1}\mathcal{H}_\alpha,& \ \alpha\in(0,2),\\
\mathcal{P}_{\alpha,a}^{ct}[0,\IF),& \ \alpha=2.
\end{array}
\RT.
\EQNY
Since for all $u$ large
 \BQNY
 \pk{u\tau_u\leq x \Big| \tau_u\leq T}=\frac{\pk{\sup_{t\in[0,u^{-1}x]}(X(t)-g(t))>u}}{\pk{\sup_{t\in[0,T]}(X(t)-g(t))>u}},
 \EQNY
then using \neprop{ruintimelocal}, we obtain for $x\in(0,\IF)$
\begin{align*}
\pk{u\tau_u\leq x \Big| \tau_u\leq T}\sim
\LT\{\begin{array}{ll}
\frac{\int_{0}^{x}e^{-ct}dt}{\int_{0}^{\IF}e^{-ct}dt} ,& \ \alpha\in(0,2),\\
\frac{\mathcal{P}_{\alpha,a}^{ct}[0,x]}{\mathcal{P}_{\alpha,a}^{ct}[0,\IF)},& \ \alpha=2.
\end{array}
\RT.
\end{align*}
%Then the result follows.\\
\COM{\underline{Case 2) $g(t)=-ct^2$:} Applying \neprop{Thm3} statement i) we obtain
\BQNY
\pk{\max_{t\in[0,T]}\LT(X(t)-ct^2\RT)>u}\sim\frac{\sqrt{\pi}}{2}a^{1/\alpha}c^{-1/2}u^{2/\alpha-1/2}\mathcal{H}_\alpha\Psi(u),\quad u\rw\IF.
\EQNY
\COM{In fact, $u+c|t-T|^2$ attains its minimum at $t_0=T$ with $u$, then
\BQNY
\lim_{u\rightarrow\IF}\sup_{t\not=T, \ |t-T|<o(1)}\left|\frac{\frac{u+c|t-T|^2}{u}-1}{\frac{c}{u}|t-T|^2}-1\right|= 0.
\EQNY
and
$$r_\eta(t)=1-a|t-T|^\alpha +o(|t-T|^\alpha), \ t\rw T.$$}
By \neprop{ruintimelocal}, we obtain for $x\in(0,\IF)$
\begin{align*}
\pk{u^{1/2}\tau_u\leq x \Big| \tau_u\leq T}\sim\frac{\int_{0}^{x}e^{-ct^2}dt}{\int_{0}^{\IF}e^{-ct^2}dt}.
\end{align*}}
{\bf Proof of Example \ref{exam3}:}
We have that  $X(t)=\frac{B_\alpha(t)}{\sqrt{Var(B_\alpha(t))}}$ is locally stationary with  correlation function
\BQNY
r_X(t,t+h)=\frac{\abs{t}^\alpha+\abs{t+h}^\alpha-\abs{h}^\alpha}{2\abs{t(t+h)}^{\alpha/2}}=1-\frac{1}{2 t^{\alpha}}|h|^\alpha +o(|h|^\alpha),\quad h\rw 0
\EQNY
for any $t>0$.  Since $g(t)=c\sin \LT(\frac{2\pi t}{T}\RT), t\in[T, (n+1)T]$ attains its maximum at $t_j=\frac{(4j+1)T}{4}, j\le n$ and
\BQNY
g(t)=c-2c\LT(\frac{\pi}{T}\RT)^2|t-t_j|^2(1+o(1)),\  t\rw t_j,\ j\le n
\EQNY
the claim follows by applying  Remarks \ref{glass} statement i). \QED

%In this part, all proofs of the examples  in section 5 are given.

{\bf Proof of Example \ref{exam1}:}
First note that the variance function of $X(t)$ is given by $\sigma^2 (t)=t(1-t)$ and correlation function is given by $r(t,s)=\frac{\sqrt{s(1-t)}}{\sqrt{t(1-s)}},0\leq s<t\leq1$.

\underline{Case 1) The proof of \eqref{ex1eq1}:}
Clearly, $m_u(t):=\frac{\sqrt{t(1-t)}}{1+ct/u}$
 attains its maximum over $[0,1]$ at the unique point  $t_u=\frac{u}{c+2u}\in(0,1)$ which converges to $t_0=\frac{1}{2}$ as $u\rightarrow\IF$, and
$m_u^*:=m_u(t_u)=\frac{1}{2\sqrt{1+c/u}}.$ Furthermore, we have
\BQN\label{MuM}
\frac{m_u^*}{m_u(t)}-1&=&\frac{u+ct}{\sqrt{t(1-t)}}\frac{\sqrt{t_u(1-t_u)}}{u+ct_u}-1
=\frac{(u+ct)\sqrt{t_u(1-t_u)}-(u+ct_u)\sqrt{t(1-t)}}{\sqrt{t(1-t)}(u+ct_u)}\nonumber\\
&=&\frac{(u+ct)^2t_u(1-t_u)-(u+ct_u)^2t(1-t)}{\sqrt{t(1-t)}(u+ct_u)[(u+ct)\sqrt{t_u(1-t_u)}+(u+ct_u)\sqrt{t(1-t)}]}.
%&=:& {I_1}/{I_2}.
\EQN
Setting $\Delta(u)=\LT[-\frac{(\ln u)^q}{u},\frac{(\ln u)^q}{u}\RT]$, and $(t_u+\Delta(u))\subset[0,\frac{1}{2}]$ for all $u$ large, we have
\BQN\label{MuM1}
(u+ct)^2t_u(1-t_u)-(u+ct_u)^2t(1-t)
&=&u^2[(t_u-t_u^2)-(t-t^2)]+2 c u t t_u (t-t_u) + c^2 t t_u (t-t_u)\nonumber\\
%&=&u^2[(t-t_u)(t+t_u)-(t-t_u)]+2 c u t t_u(t-t_u)+c^2tt_u(t-t_u)\nonumber\\
%&=&(t-t_u)[u^2((t+t_u)-1)+2 c u t t_u+c^2tt_u]\nonumber\\
%&=&(t-t_u)[u^2((t+t_u)-1)+c t u]\nonumber\\
%&=&(t-t_u)u[u(t+t_u)-c t_u-2 u t_u+ c t]\nonumber\\
%&=&(t-t_u)u[ut-ut_u+c(t-t_u)]\nonumber\\
&=&(t-t_u)^2u(u+c)
\EQN
and
\BQNY
\frac{u^4}{2\LT(u+\frac{c}{2}\RT)^2}-u^{-1/2}\leq 2(u+ct)^2[t(1-t)]\leq \frac{1}{2}\LT(u+\frac{c}{2}\RT)^2
\EQNY
for all $t\in (t_u+\Delta(u))$. Then
%%%%%%%%%%%%%%%%%%%%%%%%%%%%%%%%%%%%%%%%%%%%%%%%%%%%%%%%%%%%%%%%%%%%%%%%%%%%%%%%%%%%%%%%%%%%%%%%%%%%%%%%%5555
\COM{
\BQNY
I_2&\sim&2(u+ct_u)^2t_u(1-t_u)\\
&=&2(u+\frac{cu}{c+2u})^2\frac{u}{c+2u}\frac{c+u}{c+2u}\\
&=&2(\frac{2cu+2u^2}{c+2u})^2\frac{cu+u^2}{(c+2u)^2}\\
&=&2u^2\frac{2c+2u}{c+2u}\frac{\frac{c}{u}+1}{(\frac{c}{u}+2)^2}\\
&\rw&\frac{1}{2}u^2
\EQNY
as $u\rw\IF$.
}%%%%%%%%%%%%%%%%%%%%%%%%%%%%%%%%%%%%%%%%%%%%%%%%%%%%%%%%%%%%%%%%%%%%%%%%%%%%%%%%%%%%%%%%%%%%%%%%%%%5555
\BQN\label{Ex11}
\lim_{u\rightarrow\IF}\sup_{t\in\Delta(u), t\not=0}\left|\frac{ {m_u^*}/{m_u(t_u+t)}-1}{2t^2}-1\right|
=\lim_{u\rightarrow\IF}\sup_{t\in\Delta(u), t\not=0}\left|\frac{ {m_u^*}/{m_u(t_u+t)}-1}{2(ut)^2u^{-2}}-1\right|= 0.
\EQN
Furthermore, since
\BQNY
r(t,s)=\frac{\sqrt{s(1-t)}}{\sqrt{t(1-s)}}=1+\frac{\sqrt{s(1-t)}-\sqrt{t(1-s)}}{\sqrt{t(1-s)}}
%&=&1+\frac{s(1-t)-t(1-s)}{\sqrt{t(1-s)}(\sqrt{s(1-t)}+\sqrt{t(1-s)})}\\
=1-\frac{t-s}{\sqrt{t(1-s)}(\sqrt{s(1-t)}+\sqrt{t(1-s)})},
\EQNY
and
\BQNY
\frac{1}{2}-\frac{1}{u}\leq\sqrt{t(1-s)}(\sqrt{s(1-t)}+\sqrt{t(1-s)})\leq \frac{1}{2}+\frac{1}{u}
\EQNY
for all $s<t,\ s,t\in (t_u+\Delta(u))$,
we have
\BQNY
\lim_{u\rightarrow\IF}\underset{t\not=s}{\sup_{t,s\in\Delta(u)}}\left|\frac{1-r(t_u+t,t_u+s)}{2|t-s|}-1\right|=0.%\ \ {\bf (Assumption \ A3)}
\EQNY
Next for some small $\theta\in(0,\frac{1}{2})$, we have
$$
\E{(\overline{X}(t)-\overline{X}(s))^2}=2(1-r(t,s))\le \frac{\abs{t-s}}{(\frac{1}{2}-\theta)^2} %\ \ {\bf (Assumption \ A5)}
$$
holds for all $s,t\in[\frac{1}{2}-\theta, \frac{1}{2}+\theta]$. Moreover, by (\ref{MuM}), (\ref{MuM1}) and
\BQNY
2(u+ct)^2[t(1-t)]\leq 2\LT[u+c\LT(\frac{1}{2}+\theta\RT)\RT]^2\LT(\frac{1}{2}+\theta\RT)^2
\EQNY
for all $t\in[\frac{1}{2}-\theta, \frac{1}{2}+\theta]$, we have that for any $t\in[\frac{1}{2}-\theta, \frac{1}{2}+\theta]\setminus(t_u+\Delta(u))$
 \BQNY
 \frac{m_u^*}{m_u(t)}-1\geq\frac{(\ln u)^{2q}}{2[u+c(\frac{1}{2}+\theta)]^2(\frac{1}{2}+\theta)^2},
 \EQNY
 and further
 \BQN\label{ex112}
  \frac{m_u^*}{m_u(t)}\geq1+\mathbb{C}_1\frac{(\ln u)^q}{u^2},
  \quad t\in[\frac{1}{2}-\theta, \frac{1}{2}+\theta]\setminus(t_u+\Delta(u)).
 \EQN
Consequently, by \neprop{MainThm1}
\BQNY
\pk{\sup_{t\in[t_0-\theta, t_0+\theta]}(X(t)-ct)>u}\sim8\mathcal{H}_1u\int_{-\IF}^{\IF}e^{-8t^2}dt\Psi\LT(2\sqrt{cu+u^2}\RT)\sim e^{-2(u^2+cu)}.
\EQNY
In addition, since $\sigma_\theta:=\max_{t\in[0,1]/[t_0-\theta, t_0+\theta]}\sigma(t)<\sigma(t_0)=\frac{1}{2}$, by Borell-TIS inequality
\BQN\label{ex111}
\pk{\sup_{t\in[0,1]\setminus[t_0-\theta, t_0+\theta]}(X(t)-ct)>u}&\le& \pk{\sup_{t\in[0,1]\setminus[t_0-\theta, t_0+\theta]} X(t) > u}
\le\exp\LT(-\frac{\LT(u-\E{\sup_{t\in[0,1]}X(t)}\RT)^2}{2\sigma_\theta^2}\RT)\nonumber\\
&=&o(e^{-2(u^2+cu)}).
\EQN
Thus, by the fact that
\BQNY
&&\pk{\sup_{t\in[0,1] }(X(t)-ct)>u}\geq\pk{\sup_{t\in[t_0-\theta, t_0+\theta]}(X(t)-ct)>u}
\EQNY
and
\BQNY
\pk{\sup_{t\in[0,1] }(X(t)-ct)>u}\le\pk{\sup_{t\in[t_0-\theta, t_0+\theta]}(X(t)-ct)>u} +
\pk{\sup_{t\in[0,1]\setminus[t_0-\theta, t_0+\theta]}(X(t)-ct)>u},
\EQNY
we conclude that
$$
\pk{\sup_{t\in[0,1]}(X(t)-ct)>u}\sim e^{-2(u^2+cu)}.
$$
For any $u>0$
\BQNY
 \pk{u\LT(\tau_u-\frac{u}{c+2u}\RT)\leq x \Big| \tau_u\leq 1}=\frac{\pk{\sup_{t\in[0,t_u+u^{-1}x]}(X(t)-ct)>u}}{\pk{\sup_{t\in[0,1]}(X(t)-ct)>u}}
\EQNY
\COM{then by \eqref{ex111}
\begin{align*}
 \pk{u\LT(\tau_u-\frac{u}{c+2u}\RT)\leq x \Big| \tau_u\leq 1}\sim\frac{\pk{\sup_{t\in[t_0-\theta,t_u+u^{-1}x]}(X(t)-ct)>u}}{\pk{\sup_{t\in[t_0-\theta,t_0+\theta]}(X(t)-ct)>u}},
\end{align*}}
 and by \netheo{PreThm1}
 \BQNY
 \pk{\sup_{t\in\LT[t_u-\frac{(\ln u)^q}{u},t_u+u^{-1}x\RT]}(X(t)-ct)>u}\sim8\mathcal{H}_1u\int_{-\IF}^{x}e^{-8t^2}dt\Psi\LT(2\sqrt{cu+u^2}\RT).
 \EQNY
The above combined with \eqref{ex112} and \eqref{ex111} implies that as $u\rw\IF$
\BQNY
\pk{\sup_{t\in[0,t_u+u^{-1}x]}(X(t)-ct)>u}\sim\pk{\sup_{t\in\LT[t_u-\frac{(\ln u)^q}{u},t_u+u^{-1}x\RT]}(X(t)-ct)>u}
\sim8\mathcal{H}_1u\int_{-\IF}^{x}e^{-8t^2}dt\Psi\LT(2\sqrt{cu+u^2}\RT).
\EQNY
Consequently,
\BQNY
\pk{u\LT(\tau_u-\frac{u}{c+2u}\RT)\leq x \Big| \tau_u\leq 1}\sim\frac{\int_{-\IF}^{x}e^{-8t^2}dt}{\int_{-\IF}^{\IF}e^{-8t^2}dt}=\Phi(4x), \ \ x\in(-\IF,\IF).
\EQNY
\underline{Case 2) The proof of \eqref{ex1eq2}:} We have $t_u=\frac{u}{c+2u}\in(0,\frac{1}{2})$ which converge to $t_0=\frac{1}{2}$ as $u\rw \IF$. Since
$$\frac{1}{2}-t_u\sim\frac{c}{4u}, \quad u\rw\IF,$$
  by \neprop{MainThm1}
$$
\pk{\sup_{t\in[0,1/2]}(X(t)-ct)>u}\sim8\mathcal{H}_1u\int_{-\IF}^{c/4}e^{-8t^2}dt\Psi\LT(2\sqrt{cu+u^2}\RT)\sim \Phi(c)e^{-2(u^2+cu)}.
$$
As for the proof of \underline{Case 1)} we obtain further
 \BQNY
 \pk{u\LT(\tau_u-\frac{u}{c+2u}\RT)\leq x \Big| \tau_u\leq \frac{1}{2}}\sim\frac{\int_{-\IF}^{x}e^{-8t^2}dt}{\int_{-\IF}^{c/4}e^{-8t^2}dt}\sim\Phi(4x)/\Phi(c), \ \ x\in(-\IF,c/4].
 \EQNY
\underline{Case 3) The proof of \eqref{ex1eq3}:} We have that $\sigma(t)$ attains its maximum over $[0,1]$ at the unique point $t_0=\frac{1}{2}$, which is also the unique maximum point of $\frac{c}{2}-c\LT|t-\frac{1}{2}\RT|,t\in[0,1]$. Furthermore,
$$
\sigma (t)= \sqrt{t(1-t)}\sim\frac 1 2-\LT(t-\frac{1}{2}\RT)^2,\ \ t\to\frac{1}{2}
$$
and
$$
r(t,s)\sim1-2|t-s|,\ \ \ s,t\to  \frac{1}{2}.
$$
By \neprop{Thm2} as $u\rw\IF$
\BQNY
\pk{\sup_{t\in[0,1]}\LT(X(t)+\frac{c}{2}-c\LT|t-\frac{1}{2}\RT|\RT)>u}\sim
8\mathcal{H}_1u\int_{-\IF}^{\IF}e^{-\LT(8|t|^2+4c|t|\RT)}dt\Psi\LT(2u-c\RT)
\sim  2\Psi(c) e^{-2(u^2-cu)}
\EQNY
and in view of  \neprop{ruintime} ii)
\begin{align*}
\pk{u\LT(\tau_u-\frac{1}{2}\RT)\leq x \Big| \tau_u\leq 1}\sim\frac{\int_{-\IF}^{x}
e^{-\LT(8|t|^2+4c|t|\RT)}dt}{\int_{-\IF}^{\IF}e^{-\LT(8|t|^2+4c|t|\RT)}dt},\ u\rw\IF.
\end{align*}

 \QED

\COM{\BEL\label{lem1}
Let $f\in C_0^*(E)$ and $E$ is a closed interval in $\R$, then $\int_{E} e^{-cf(t)}dt<\IF$ holds for any constant $c>0$.
\EEL
\prooflem{lem1} Since by $f\in C_0^*(E)$ the function $f$ is locally bounded and integrable, it suffice to show that
$\int_b^\IF e^{-cf(t)}dt<\IF$ and $\int_{-\IF}^{-b} e^{-cf(t)}dt<\IF$  for some $b$ large enough. By the assumption, we can find $a\in(0,\IF)$, such that
$cf(t)>t^{\epsilon_1}$ holds for any $t>a$. Hence
\BQNY
\int_a^\IF e^{-c_1f(t)}dt\leq
\int_{0}^\IF e^{-t^{\epsilon_1}}dt=\frac{1}{\epsilon_1}\Gamma\LT(\frac{1}{\epsilon_1}\RT)
<\IF.
\EQNY
Similarly, we have $\int_{-\IF}^{-b} e^{-cf(t)}dt<\IF$.
Thus we establish the claim. \QED}
\COM{
Let $\{X(t_1,t_2),(t_1,t_2)\in[0,T_1]\times [0,T_2]\}$ be a centered Gaussian field with continuous trajectories,  standard deviation function $\sigma(\cdot,\cdot)$ and correlation function $r(\cdot,\cdot,\cdot,\cdot)$.
Further, $\{g(t_1,t_2),(t_1,t_2)\in[0,T_1]\times [0,T_2]\}$ is a bounded measurable function,  and we set
$$M_u(t_1,t_2):=\frac{u-g(t_1,t_2)}{\sigma(t_1,t_2)}.$$
\textbf{D1}: For any fixed large $u$, $M_u(t_1,t_2), (t_1,t_2)\in[0,T_1]\times[0,T_2]$ attains its minimum, denoted by $M_u^*$, at a unique point $(t_1(u),t_2(u))$, and
$$
\lim_{u\to\IF} (t_1(u),t_2(u))=(t_1^*,t_2^*)\in[0,T_1]\times[0,T_2].
$$
\textbf{D2}: There exist some interval $\Delta_1(u)\times\Delta_2(u)$ and some constant $\lambda_1,\lambda_2>0$ such that
$$
\lim_{u\rightarrow\IF}\underset{(t_1,t_2)\neq(0,0)}{\sup_{(t_1,t_2)\in\Delta_1(u)\times\Delta_2(u)}}\left|\frac{ M_u(t_1(u)+t_1,t_2(u)+t_2)/M_u^*-1}{[f_1(u^{\lambda_1}|t_1|)+f_2(u^{\lambda_2}|t_2|)]u^{-2}}-1\right|= 0,
$$
where  $f_1(\cdot),f_2(\cdot)$ is a positive function over $(0,\IF)$.\\ %|t|)\in\mathcal{R}_\beta,\beta>0$.
%In our analysis the following expansion of the correlation function $r(\cdot,\cdot)$ around $t_u$ is also important.
\textbf{D3}: There exist some constants $\alpha_i\in(0,2],\ i=1,2$
$$
\lim_{u\rightarrow\IF}\underset{(t_1,t_2)\neq(s_1,s_2),(t_1,t_2)\in\Delta_1(u)\times\Delta_2(u)} {\sup_{(s_1,s_2)\in\Delta_1(u)\times\Delta_2(u)}}\left|\frac{1-r(t_1(u)+t_1,t_1(u)+s_1,t_2(u)+t_2,t_2(u)+s_2)}{\rho_1^2(|t_1-s_1|)+\rho_2^2(|t_2-s_2|)}-1\right|=0,
$$
with $\Delta_1(u)\times\Delta_2(u)$ given in \textbf{A2} and $\rho_i^2\in\mathcal{R}_{\alpha_i}, i=1,2$.\\
\textbf{D4}: For any $(t_1,t_2)\in [0,T_1]\times[0,T_2]\setminus((t_1(u),t_2(u))+\Delta^1_u\times\Delta^2_u)$
$$M_u(t_1,t_2)\geq M_u(s_1,s_2)$$
holds for any $(s_1,s_2)\in(t_1(u),t_2(u))+\Delta^1_u\times\Delta^2_u$ when $u$ is large enough.\\
\textbf{D5}: For some positive $G, \varsigma_1,\varsigma_2$,
$$\E{(\overline{X}(t_1,t_2)-\overline{X}(s_1,s_2))^2}\leq G(|t_1-s_1|^{\varsigma_1}+|t_2-s_2|^{\varsigma_2})$$
holds for all $(s_1,s_2),(t_1,t_2)\in[0,T_1]\times[0,T_2]$.
\BT
For $T_1,T_2>0$, let $\{X(t_1,t_2),(t_1,t_2)\in[0,T_1]\times [0,T_2]\}$ be the centered Gaussian field  and $g(t_1,t_2)$ is a bounded measurable function over $[0,T_1]\times [0,T_2]$.
If \textbf{D1}-\textbf{D5} are satisfied and $f_1,\ f_2$ satisfies $\textbf{B1}$. Set $\lim_{t\downarrow 0}\frac{\rho_i^2(t)}{t^{2/{\lambda_i}}}=\eta_i\in[0,\IF]$, $h_i(t)=\frac{1}{\sigma_0^2}f_i(\sigma_0^{2/{\alpha_i}}\eta_i^{-1/\alpha} |t|) $ and $\sigma_0=\sigma(t_1^*,t_2^*)$. Then for $X_g(t_1,t_2) =X(t_1,t_2)+ g(t_1,t_2)$, we have
\BQN
\pk{\sup_{(t_1,t_2)\in[0,T_1]\times[0,T_2]}X_g(t_1,t_2)>u}\sim  \prod_{i=1}^{2}H_i(\alpha_i,\lambda_i,d_i^*,f_i,\rho_i,\eta_i,u,\sigma_0)\Psi(M_u^*),
\EQN
where $M_u^*:=M_u(t_1^*,t_2^*)$,
\BQNY
H_i(\alpha_i,\lambda_i,d_i^*,f_i,\rho_i,\eta_i,u,\sigma_0)=\left\{
\begin{array}{ll}
\frac{\sigma_0^{-2/\alpha}\mathcal{H}_{\alpha_i}}{u^{\lambda_i}\orh_i(u^{-1})}\int_{-d_i^*}^{\IF}e^{-\frac{f_i(|t|)}{\sigma^2_0}}dt,&\hbox{if} \ \ \eta_i=\IF,\\
\mathcal{P}_{\alpha_i}^{h_i}[-\sigma_0^{-2/{\alpha_i}}{\eta_i}^{1/{\alpha_i}}d_i^*,\IF],& \hbox{if} \ \ \eta_i\in(0,\IF),\\
1,&  \hbox{if} \ \ \eta_i=0,
\end{array}
 \right.\
d_i^*:=
\left\{
\begin{array}{ll}
0, & \text{if  \textbf{C1'} holds},\\
d_i, & \text{if \textbf{C2'} holds},\\
\IF, & \text{if \textbf{C3'} holds},
\end{array}
\right.
 \EQNY
and \textbf{C1'}: $t_i(u)\sim d_iu^{-\nu_i}$, or $T_i-t_i(u)\sim d_iu^{-\nu_i}$, $\nu_i>\lambda_i$;\\
\textbf{C2'}: $t_i(u)\sim d_iu^{-\nu_i}$, or $T_i-t_i(u)\sim d_iu^{-\nu_i}$, $\nu_i=\lambda_i$;\\
\textbf{C3'}: $t_i(u)\sim d_iu^{-\nu_i}$, or $T_i-t_i(u)\sim d_iu^{-\nu_i}$, $\nu_i<\lambda_i$ or $t_i^*\in(0,T_i)$,\\
$\Delta_i(u)=[\delta_i(u),\delta^i_u]$, $\delta_i(u)=d_iu^{-\nu_i}$ if \textbf{C1'} and \textbf{C2'} hold,
$\delta_i(u)=\delta^i_u=\LT(\frac{(\ln u)^{q_i}}{u}\RT)^{\lambda_i}$ for some large positive constant constant $q_i$ if \textbf{C3'} holds.
\ET
}
\section*{Acknowledgments} % We are thankful to the referees for several suggestions which improved our manuscript.
Thanks to Swiss National Science Foundation Grant no. 200021-166274.
KD acknowledges partial support by NCN Grant No 2015/17/B/ST1/01102 (2016-2019).

{
\footnotesize
\bibliographystyle{ieeetr}
\bibliography{trendPaper}
}
\end{document}